\documentclass[11pt]{article}
\usepackage{amsmath,verbatim,amssymb,array}
\usepackage[cp1250]{inputenc}
\usepackage[T1]{fontenc}
\usepackage[english]{babel}
\usepackage{a4wide,cite}
\usepackage{graphicx}
\usepackage{bm}
\usepackage{color}

\newenvironment{pf}{{\it Proof.\ }\rm }{ \hfill$\Box$\ \\[0.25ex]}
\numberwithin{equation}{section}
\newtheorem{thm}{Theorem}[section]
\newtheorem{lem}[thm]{Lemma}

\newtheorem{alg}[thm]{Algorithm}
\newtheorem{prob}[thm]{Problem}

\newenvironment{rem}{\refstepcounter{thm}\ \\[2mm]%
                  \noindent{\it Remark~\thethm\ \ }}{\ \\[2mm]\rm}

\newcommand{\intT}{\mbox{$\displaystyle\int\!\!\!\!\int\limits_{\!\!\!\!T}$}}
\newcommand{\N}{\mathbb N}

\newcommand{\R}{\mathbb R}
\newcommand{\Z}{\mathbb Z}
\newcommand{\Zplus}{\Z_+}
\newcommand{\C}{\mathbb C}

\newcommand{\bl}[1]{\mbox{$\bm{#1}$}}
\newcommand{\sbl}[1]{\mbox{\scriptsize$\bm{#1}$}}  
\newcommand{\tbl}[1]{\mbox{\tiny$\bm{#1}$}}   
\newcommand{\bla}{\mbox{$\bm{\alpha}$}}
\newcommand{\lbla}{\mbox{$\left|\bla\right|$}}
\newcommand{\sbla}{\mbox{\scriptsize$\bla$}}

\newcommand{\blc}{\mbox{$\bl c$}}
\newcommand{\sblc}{\mbox{$\sbl c$}}
\newcommand{\lblc}{\mbox{$\left|\blc\right|$}}

\newcommand{\dx}{\,\mbox{{\rm d}$\bl x$}}
\newcommand{\blk}{\mbox{$\bm{k}$}}
\newcommand{\sblk}{\mbox{\scriptsize$\blk$}}
\newcommand{\bll}{\mbox{$\bm{\bm l}$}}
\newcommand{\sbll}{\mbox{\scriptsize$\bll$}}
\newcommand{\blm}{\mbox{$\bm{\mu}$}}

\newcommand{\sP}{\mbox{$\mathsf P$}}
\newcommand{\sR}{\mbox{$\mathsf R$}}
\newcommand{\sQ}{\mbox{$\mathsf Q$}}
\newcommand{\sS}{\mbox{$\mathsf S$}}
\newcommand{\sT}{\mbox{$\mathsf T$}}
\newcommand{\sW}{\mbox{$\mathsf W$}}

\newcommand{\hts}{\hspace*{0.1ex}}
\newcommand{\BB}{\mathcal{B}}
\newcommand{\MM}{\mathcal{M}}

\newcommand{\rbinom}[2]{\mbox{$\displaystyle\binom{#1}{#2}^{\!\!-1}\!\!$}}
\newcommand{\hyper}[5]{\,{}_{#1}F_{#2}\!\left(%
            		    \begin{array}{cc}{#3}\\{#4} \end{array}\big|\,{#5} \right)}

\begin{document}
	
 \thispagestyle{empty}   

 \title{\Large Constrained  approximation of rational triangular B\'ezier surfaces\\
                   by polynomial triangular  B\'ezier surfaces}
                   
 \author{Stanis{\l}aw Lewanowicz${}^{a,\!}$%
           	 \footnote{Corresponding author\newline
		\hspace*{1.5em}\textit{Email addresses}: Stanislaw.Lewanowicz@ii.uni.wroc.pl 
		(Stanisław Lewanowicz),\newline
		Pawel.Keller@mini.pw.edu.pl (Paweł Keller), Pawel.Wozny@ii.uni.wroc.pl (Paweł Woźny)}
		\,\:, Pawe{\l} Keller${}^{b}$, Pawe{\l} Wo\'zny${}^{a}$
		}

 \date{\small\it ${}^{a}\!\!\!$ Institute of Computer Science, University of Wroc{\l}aw,
  	  	ul.~Joliot-Curie 15, 50-383 Wroc{\l}aw, Poland\\	        
 \noindent ${}^{b}\!\!\!$ Faculty of Mathematics and Information Science, 
 		Warsaw University of Technology,\\ ul. Koszykowa 75, 00-662 Warszawa, Poland\\[2ex]}
   
\maketitle

\noindent {\small\textbf{Abstract}.
 		  We propose a novel approach to the problem of polynomial approximation
of rational B\'ezier triangular patches with prescribed boundary control points.
The method is very efficient thanks to using recursive properties
of the bivariate dual Bernstein polynomials and applying a~smart algorithm
for evaluating a collection of two-dimensional integrals.  
		 Some illustrative examples are given.}\\
	 


\section{Introduction and preliminaries}                        \label{S:intro} 

Rational triangular B\'ezier surfaces are an important tool in computer graphics. 
However, they may be sometimes inconvenient  in practical applications.
The reason is that evaluation of integrals or derivatives of rational expressions is cumbersome. Also, 
it happens that a rational surface produced in one CAD system is to be imported into another system 
which can handle only polynomial surfaces.

In order to solve the two problems above, different algorithms for approximating the rational surface by 
polynomial surface are proposed \cite{CW11,Hu13,Sha13,XW09,XW10,ZW06}. The spectrum of methods
contains hybrid algorithm \cite{ZW06}, progressive iteration approximation
\cite{CW11,Hu13}, least squares approximation and linear programming \cite{Sha13},
and approximation by B\'ezier surfaces with control points obtained by successive degree elevation
of the rational B\'ezier surface \cite{XW09,XW10}.  As a rule, no geometric 
constraints are imposed, which means a serious drawback: if we start with a patchwork of smoothly
connected rational B\'ezier triangles and approximate each patch separately, we do not obtain a smooth composite surface.

In this paper, we propose a method for solving the problem of the constrained least squares approximation 
of a rational triangular B\'ezier patch by a polynomial triangular B\'ezier patch; see Problem~\ref{P:main} below. 
The method is based on the idea of using constrained dual bivariate Bernstein polynomials.
Using a fast recursive scheme of evaluation of B\'ezier form coefficients of the bivariate dual Bernstein polynomials, and applying a swift adaptive scheme of numerical computation of a collection 
of double integrals involving rational functions resulted in high efficiency of the method.

The outline of the paper is as follows. Section~\ref{S:main} brings a complete
solution to the approximation problem. 
Some comments on the algorithmic implementation of the method 
are given in Section~\ref{S:impl}; some technical details of the implementation are presented in Appendix A.
In Section~\ref{S:exmp}, several examples are given to show the efficiency of the method. 
In Appendix~B, some basic information on the Hahn orthogonal polynomials is recalled.
	
We end this section by introducing some notation.
For  $\bl y:=(y_1,y_2,\ldots,y_d)\in\R^d$, we  denote 
$|\bl y|:=y_1+y_2+\ldots+y_d$,  and $\|\bl y\|:=\left(y^2_1+y^2_2+\ldots+y^2_d\right)^{\frac12}$.

\begin{figure}[htb]
	\begin{center}     
        \begin{tabular}{ccc}
	\hspace*{-1cm} 
	\includegraphics[scale=0.225,angle=-90]{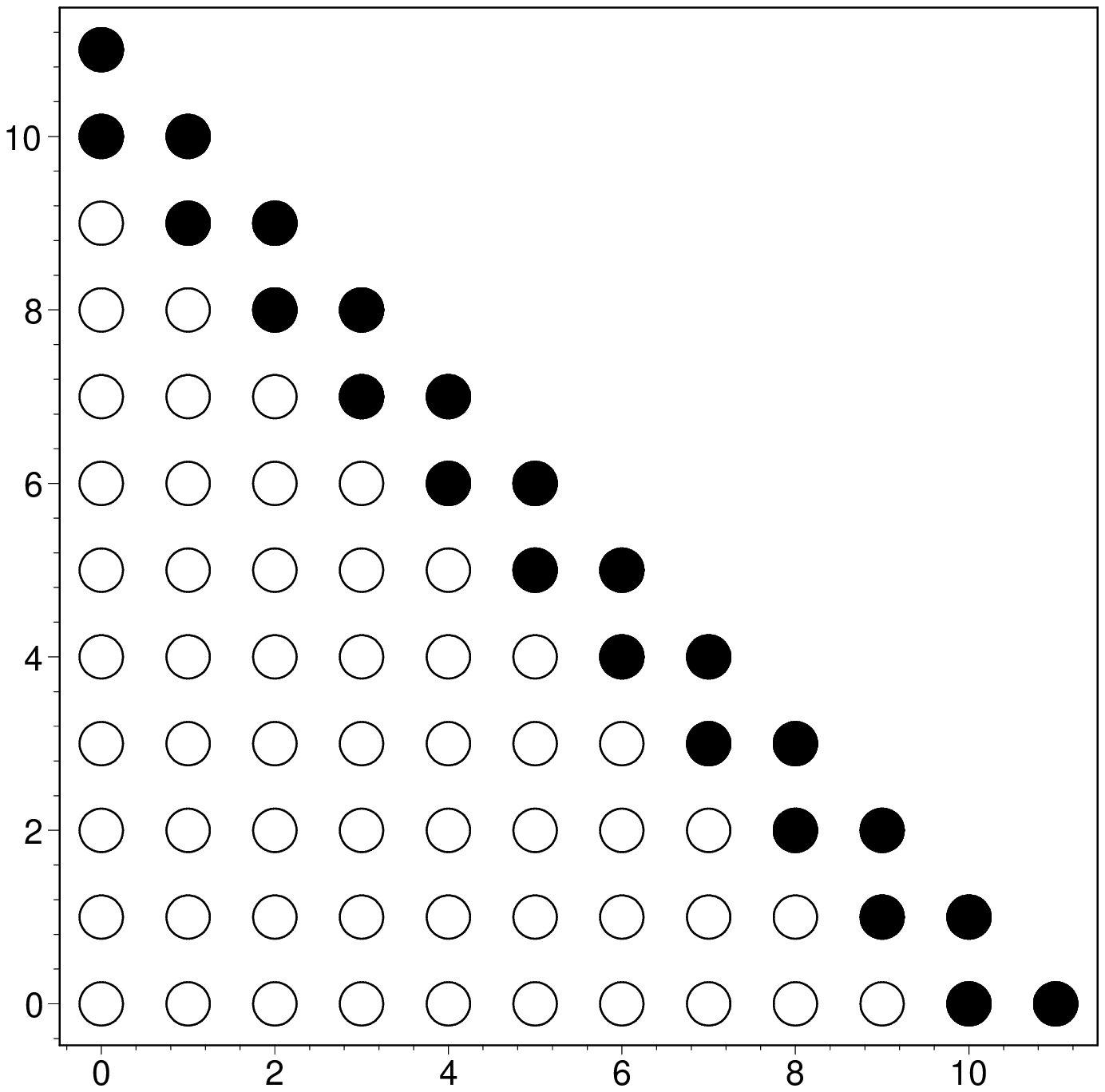} &
	\hspace*{-1.8cm} 
	\includegraphics[scale=0.225,angle=-90]{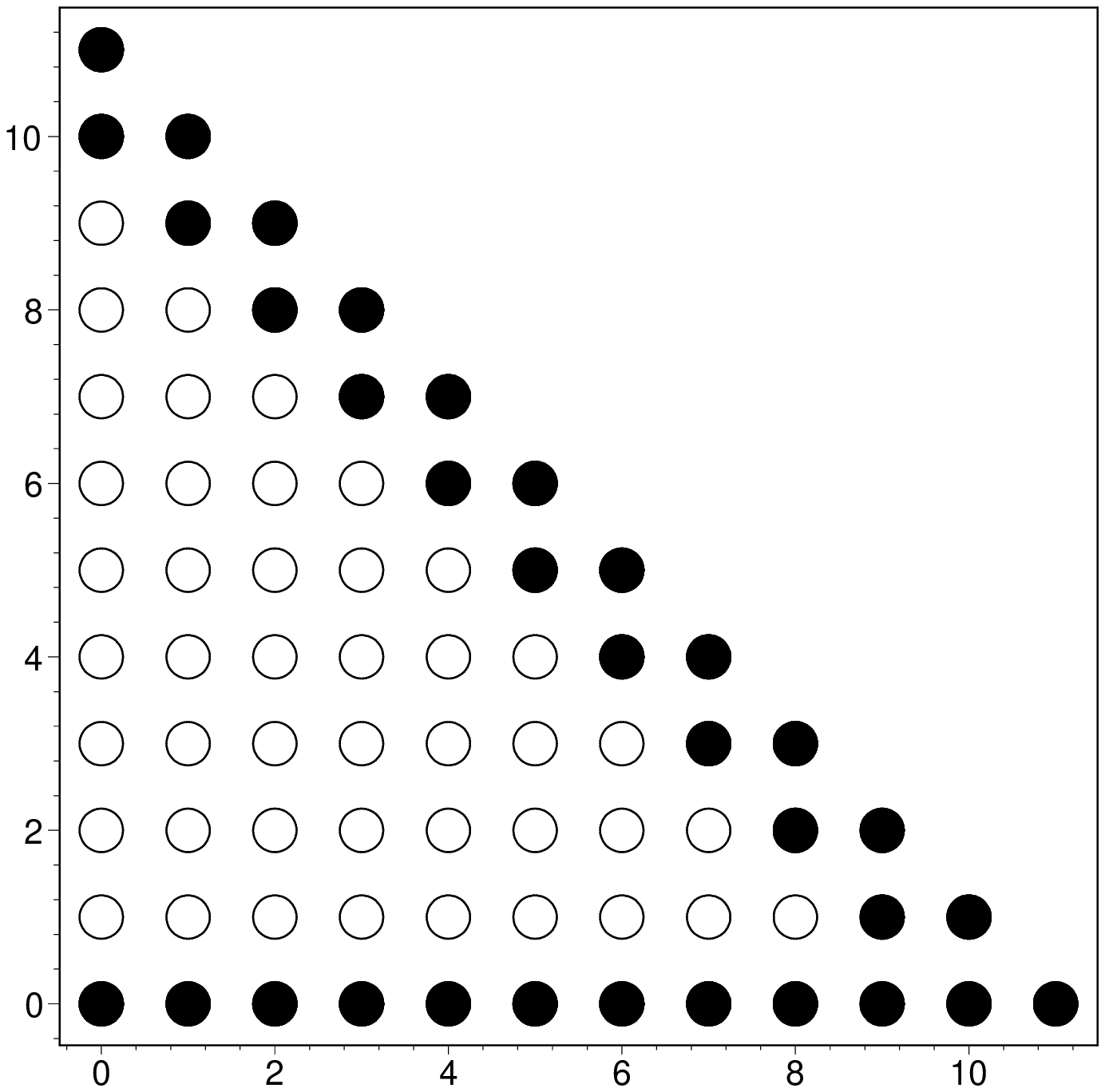}&
	\hspace*{-1.8cm} 
	\includegraphics[scale=0.225,angle=-90]{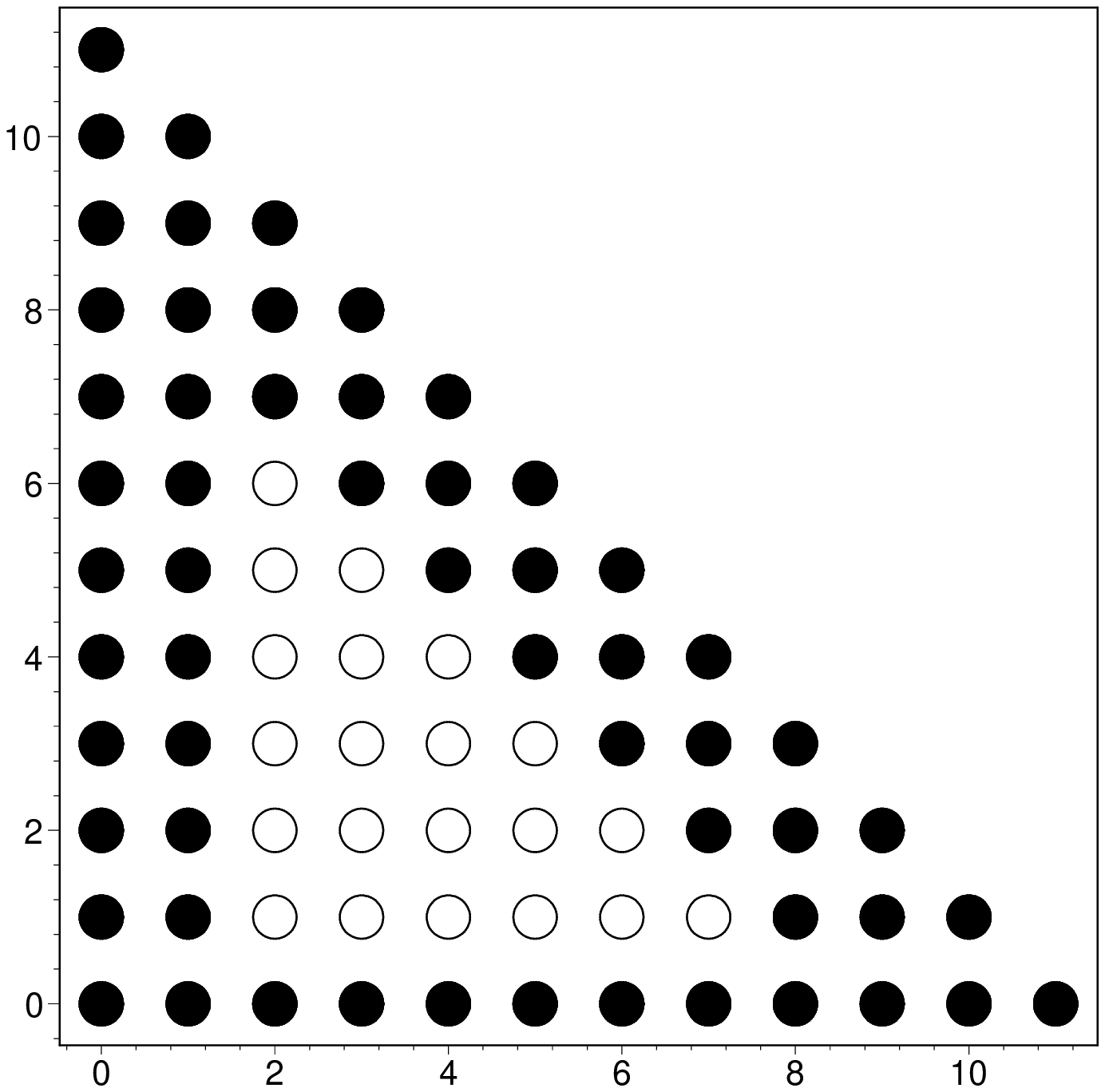}
	\\[-0.2ex]
	\hspace*{-0.9cm} $\bm{c}=(0,0,2)$ & \hspace*{-1.65cm} $\bm{c}=(0,1,2)$ 
	& \hspace*{-1.65cm}	$\bm{c}=(2,1,3)$ 
	\end{tabular}	 		
	\vspace*{0.125cm}
	\caption{\small Examples of sets \eqref{E:TOG} ($n=11$). Points of the set
		$\Omega^{\sbl c}_n$ 
	are marked by white discs, while the points of the set
	$\Gamma^{\sbl c}_n$ -- by black discs. Obviously, 
	$\Theta_n=\Omega^{\sbl c}_n\cup\Gamma^{\sbl c}_n$.%
	\label{fig:Fig-bc}}
	\end{center}
	\end{figure}

For $n\in\N$ and  $\bl c:=(c_1,\,c_2,\,c_3)\in\N^3$ such that $\lblc<n$, we define 
the following sets (cf. Figure~\ref{fig:Fig-bc}):
\begin{equation}\label{E:TOG}
	\left.\begin{array}{l}
	\Theta_n:=\{\bl k=(k_1,k_2)\in\N^2:\: 0\le|\bl k|\le n\},  \\[1ex]
	\Omega^{\sbl c}_n:=\{\bl k=(k_1,k_2)\in\N^2:\:k_1\ge c_1,\,k_2\ge c_2,\,|\bl k|\le n-c_3\},\\[1ex]
	\Gamma^{\sbl c}_n:=\Theta_n\setminus\Omega^{\sbl c}_n.
	  \end{array}\;\right\}
\end{equation}	

Throughout this paper, the symbol $\Pi^2_n$  denotes the space
 of all polynomials of two variables,  of total degree at most $n$.

Let $T$ be the standard triangle in $\R^2$,
\begin{equation}
	\label{E:T}
		T:=\{(x_1,\,x_2)\,:\,x_1,\,x_2\ge0,\: x_1+x_2\le1\}.		
\end{equation}	
	For $n\in\N$, and $\bl k:=(k_1,k_2)\in\Theta_n$, we denote,
\[
	\binom{n}{\bl k}:=\frac{n!}{k_1!k_2!(n-|\bl k|)!}.
\]
The \textit{shifted factorial} is defined for any $a\in\C$ by
\[
	(a)_0:=1;\qquad (a)_k:=a(a+1)\cdots(a+k-1), \qquad  k\ge1.
\]		
The \textit{Bernstein polynomial basis} in $\Pi^2_n$, $n\in\N$, is given by
 (see, e.g., \cite{Far86}, or \cite[\S 17.3]{Far02}),
\begin{equation}\label{E:Ber2}
	B^n_{\sbl k}(\bl x):=\binom{n}{\bl k}x_1^{k_1}x_2^{k_2}(1-|\bl x|)^{n-|\sbl k|},
	\qquad \bl k:=(k_1,k_2)\in\Theta_n,\quad\bl x:=(x_1,x_2).
\end{equation}

The (unconstrained) \textit{bivariate dual Bernstein basis polynomials} \cite{LW06}, 
\begin{equation}\label{E:dBer2}
	D^n_{\sbl k}(\bl\cdot;\bla)\in\Pi^2_n,	
	\qquad \bl k\in\Theta_n,
\end{equation}
are defined so that
\[  
	\left\langle D^{n}_{\sbl k},\,B^n_{\sbl l} \right\rangle_{\sbla}=\delta_{\sbl k, \sbl l},
	\qquad \bl k,\bl l\in\Theta_n.
\] 
Here $\delta_{\sbl k, \sbl l}$ equals 1 if $\bl k= \bl l$,  and 0 otherwise, while the inner product is defined by
\begin{equation}\label{E:Jinprod}
\langle f, g \rangle_{\sbla} 
:=\intT w_{\sbla}(\bl x)f(\bl x)\,g(\bl x)\dx,
\end{equation} 
where the weight function $w_{\sbla}$ is given by
\begin{equation}
	\label{E:w}
	w_{\sbla}(\bl x):=A_{\sbla}x_1^{\alpha_1}x_2^{\alpha_2}(1-|\bl x|)^{\alpha_3}, 
	\qquad\bla:=(\alpha_1,\,\alpha_2,\,\alpha_3), \quad\alpha_i>-1,
		\end{equation}	
	with \(  
	A_{\sbla}:=\Gamma(\lbla+3)/[\Gamma(\alpha_1+1)\Gamma(\alpha_2+1)\Gamma(\alpha_3+1)].
\)

For $n\in\N$ and  $\bl c:=(c_1,\,c_2,\,c_3)\in\N^3$ such that $\lblc<n$, define 
 the constrained bivariate  polynomial space 
\[  
	\Pi^{\sbl c,\,2}_n :=\left\{P\in\Pi^2_n\::\: 
P(\bl x)= x_1^{c_1}x_2^{c_2}(1-|\bl x|)^{c_3}\cdot\,Q(\bl x),
\;\mbox{where}\;Q\in\Pi^2_{n-|\sbl c|}\right\}.
\] 
It can be easily seen that the constrained set 
$\{B^n_{\sbl k}\}_{\sbl k\in\Omega^{\sbl c}_n}$
of degree $n$ bivariate Bernstein polynomials forms a basis in this space. 
We define \textit{constrained  dual bivariate Bernstein basis polynomials}, 
\begin{equation}\label{E:constrdBer2}
	D^{(n,\sbl c)}_{\sbl k}(\bl{\cdot};\bla)\in\Pi^{\sbl c,\,2}_n,
	\qquad \bl k\in\ \Omega^{\sbl c}_n,
\end{equation}
 so that
\begin{equation}
	\label{E:dualprop}
	\left\langle D^{(n,\sbl c)}_{\sbl k},\,B^n_{\sbl l}
	\right\rangle_{\sbla}=\delta_{\sbl k, \sbl l} \quad \mbox{for}
	\quad \bl k,\:\bl l\in\ \Omega^{\sblc}_n,	
\end{equation}
where the notation of \eqref{E:Jinprod} is used. For $\bl c=(0,0,0)$, basis \eqref{E:constrdBer2}
reduces to the unconstrained basis \eqref{E:dBer2} in $\Pi^2_n$. 
Notice that the solution of the least squares approximation problem in the space $\Pi^{(\sbl c,2)}_n$ can be given in terms of the polynomials $D^{(n,\sbl c)}_{\sbl k}$. Namely, we have the following result.
\begin{lem}\label{L:bestpol}
    	    Let  $F$ be a function defined
    	    on the standard triangle $T$ (cf. \eqref{E:T}).
    	    The polynomial $S_n\in\Pi^{(\sbl c,2)}_n$,  which gives the minimum
value of the norm
\[
	\|F-S_n\|_{L_2}:=\left\langle F-S_n,F-S_n\right\rangle_{\sbla}^{\frac12},
\]
is given by
\begin{equation}
	\label{E:bestpol}
S_n=
\sum_{\sbl k\in\Omega^{\tbl c}_n}\left\langle F,D^{(n,\sbl c)}_{\sbl k}\right\rangle_{\!\sbla}\,B^n_{\sbl k}.
\end{equation}
\end{lem}
\begin{pf}
	Obviously, $S_n$ has the following representation in the  Bernstein basis of the space
	$\Pi_n^{(\sbl c,2)}$:
	\[
		S_n=\sum_{\sbl k\in\Omega^{\tbl c}_n}
		\left\langle S_n,D^{(n,\sbl c)}_{\sbl k}\right\rangle_{\!\sbla}\,B^n_{\sbl k}.		
	\]
	On the other hand, a classical characterization of the best approximation polynomial $S_n$ is that 
	$\langle {F-S_n}, Q\rangle_{\sbla}=0$ holds for
	any polynomial  $Q\in\Pi_n^{(\sbl c,2)}$ (see, e.g. \cite[Thm 4.5.22]{DB08}). 
	In particular, for $Q=D^{(n,\sbl c)}_{\sbl k}$, we obtain
\[
	\left\langle {F},{D^{(n,\sbl c)}_{\sbl k}}\right\rangle_{\!\sbla}
	=\left\langle S_n,D^{(n,\sbl c)}_{\sbl k}\right\rangle_{\!\sbla} ,\qquad \bl k\in\Omega^{\sbl c}_n.
\]
Hence, the formula \eqref{E:bestpol} follows.
\end{pf}

The coefficients  $E^{\sbl k}_{\sbl l}(\bla,\blc,n) $ in the B\'ezier form of the  dual Bernstein polynomials, 
\begin{equation}
	\label{E:Dc2inB2}
	D^{(n,\sblc)}_{\sbl k}(\bl x;\bla)=\,
	\sum_{\sbl l\in\Omega^{\sblc}_n}\,E^{\sbl k}_{\sbl l}(\bla,\blc,n) \,
	B^n_{\sbl l}(\bl x),\qquad \bl k\in\Omega^{\sblc}_n,
	\end{equation}
play important role in the proposed method. Using the duality property \eqref{E:dualprop},
we obtain the following expression for the coefficients of the above expansion: 
\begin{equation}
	\label{E:Ek}
	E^{\sbl k}_{\sbl l}(\bla,\blc,n)  =\left\langle D^{(n,\sbl c)}_{\sbl k}, 
	           D^{(n,\sbl c)}_{\sbl l}\right\rangle_{\sbla}.
	\end{equation}
In a recent paper \cite{LKW15}, an efficient  algorithm was obtained
for evaluating all these coefficients  for $\blk,\,\bll\in\Omega^{\sblc}_n$, with the computational complexity $O(n^4)$, i.e., proportional to the total number of these coefficients. See Section~\ref{SS:impl-E} for details.

\section{Polynomial approximation of B\'ezier triangular surfaces with constraints}
						\label{S:main}

In this paper, we consider the following approximation problem. 
\begin{prob} \label{P:main}
	Let $\sR_n$ be a rational triangular B\'ezier surface of degree $n$,
\[  
	\sR_n(\bl x)	:=\frac{\sQ_n(\bl x)}{\omega(\bl x)}
	=\frac{\displaystyle \sum_{\sbl k\in\Theta_n}\omega_{\sbl k}r_{\sbl k}B^n_{\sbl k}(\bl x)}
	{\displaystyle \sum_{\sbl k\in\Theta_n}\omega_{\sbl k}B^n_{\sbl k}(\bl x)},\qquad \bl x\in T,
\] 
 with the control points $r_{\sbl k}\in\R^3$ 
and positive weights $\omega_{\sbl k}\in\R$, $\bl k\in\Theta_n$.
Find a~B\'ezier triangular surface of degree $m$,
\[  
\sP_m(\bl x)
:= \sum_{\sbl k\in\Theta_m}p_{\sbl k}B^m_{\sbl k}(\bl x),\qquad \bl x\in T,	         
\]  
	with the control points $p_{\sbl k}\in\R^3$, 	satisfying the conditions
\begin{equation}\label{E:gcond}
	p_{\sbl k}=g_{\sbl k}\quad \mbox{for}\quad \bl k\in\Gamma^{\sbl c}_m,
\end{equation}	
$g_{\sbl k}\in\R^3$ being prescribed control points, and  $\bl c:=(c_1,c_2,c_3)\in\N^3$ being a given
parameter vector with $|\bl c|<m$,
such that the distance between the surfaces $\sR_n$ and $\sP_m$,
\begin{equation}
	\label{E:dist}
	d(\sR_n,\sP_m):=\intT w_{\sbla}(\bl x)\|\sR_n(\bl x)-\sP_m(\bl x)\|^2\dx,
	\end{equation}
reaches the minimum.
\end{prob}

\begin{rem}\label{R:gcond}
Remember that continuity conditions for any two adjacent triangular  B\'ezier patches are given in terms 
of several rows of the control net "parallel" to the  control polygon of their common boundary
(see, e.g., \cite[Section 17]{Far02}). Therefore, constraints  \eqref{E:gcond} are natural,
in a sense (cf.  Fig.~\ref{fig:Fig-bc}). In Section~\ref{S:exmp}, we give several examples
of practical usage of this approach.
\end{rem} 

Clearly, the B\'ezier triangular patch being the solution of Problem~\ref{P:main} can be obtained in a componentwise
way. Hence it is sufficient to give a method for solving the above problem
in the case where $\sR_n$ and $\sP_m$ are scalar functions, and $g_{\sbl k}$ are numbers.

All the details of the proposed method are given in the following theorem. 

\begin{thm}
	\label{T:main}
Given the coefficients $r_{\sbl k}$ 
and positive weights $\omega_{\sbl k}$, $\bl k\in\Theta_n$, of the rational 
function
\begin{equation}
	\label{E:Rsc}
	\sR_n(\bl x):=\frac{\sQ_n(\bl x)}{\omega(\bl x)}
	=\frac{\displaystyle \sum_{\sbl k\in\Theta_n}\omega_{\sbl k}r_{\sbl k}B^n_{\sbl k}(\bl x)}
	          {\displaystyle \sum_{\sbl k\in\Theta_n}\omega_{\sbl k}B^n_{\sbl k}(\bl x)},
\end{equation}
the coefficients $p_{\sbl k}$ of the degree $m$ polynomial	
\begin{equation}\label{E:Psc}
	\sP_m(\bl x):= \sum_{\sbl k\in\Theta_m}p_{\sbl k}B^m_{\sbl k}(\bl x),
\end{equation}
minimising the error
\begin{equation}\label{E:er-sc}
	\|\sR_n-\sP_m\|^2_{L_2}:=\langle \sR_n-\sP_m,\sR_n-\sP_m\rangle_{\sbl \alpha},	
\end{equation}
with the constraints
\begin{equation}\label{E:gsc}
p_{\sbl k}=g_{\sbl k}\quad \mbox{for}\quad \bl k\in\Gamma^{\sbl c}_m,	
\end{equation}
are  given by
\begin{equation}\label{E:psc}
p_{\sbl k}=\sum_{\sbl l\in\Omega^{\sbl c}_m}\binom{m}{\bl l}\,E^{\sbl k}_{\sbl l}(\bla,\blc,m)
	              \big(u_{\sbl l}-v_{\sbl l}\big), \qquad \bl k\in\Omega^{\sbl c}_m,	
\end{equation}
where
\begin{align*}  
	u_{\sbl l}:=& \sum_{\sbl h\in\Theta_n} \binom{n}{\bl h}
	          \rbinom{n+m}{\bl h+\bl l}\,\omega_{\sbl h}r_{\sbl h}\,I_{\sbl h+\sbl l},\\ 
        \label{E:vl}
	v_{\sbl l}:=& \frac{1}{(|\bla|+3)_{2m}}\sum_{\sbl h\in\Gamma^{\sbl c}_m}\binom{m}{\bl h}	     	       
	       \left(\prod_{i=1}^{3}(\alpha_i+1)_{h_i+l_i}\right)g_{\sbl h},             	
\end{align*}
with $h_3:=m-|\bl h|$, $l_3:=m-|\bl l|$, and
\begin{equation}
	\label{E:I}
	I_{\sbl j}:=  \intT w_{\sbla}(\bl x)\frac{B^{n+m}_{\sbl j}(\bl x)}{\omega(\bl x)}\dx,
	\qquad\bl j\in\Omega^{\sblc}_{n+m}.
\end{equation}          
The symbol $E^{\sbl k}_{\sbl l}(\bla,\blc,m)$ has the meaning given in \eqref{E:Dc2inB2}.
\end{thm}

\begin{pf}
Observe that
\[
	\|\sR_n-\sP_m\|^2_{L_2}=\|\sW-\sS_m\|^2_{L_2}	
\]
where 
\[
	\sW:=\sR_n-\sT_m,\quad \sT_m:=\sum_{\sbl k\in\Gamma^{\sbl c}_m}g_{\sbl k}B^m_{\sbl k},
	\quad \sS_m:=\sum_{\sbl k\in\Omega^{\sbl c}_m}p_{\sbl k}B^m_{\sbl k},
\]
the notation being that of \eqref{E:TOG}. Thus, we want $\sS_m$ to be the best approximation polynomial
for the function $\sW$ in the space $\Pi^{\sbl c,2}_m$. Its B\'ezier coefficients are given by
\[
	p_{\sbl k}=\left\langle \sW, D^{(m,\sbl c)}_{\sbl k} \right\rangle_{\sbla}
=\sum_{\sbl l\in\Omega^{\sbl c}_m}\,E^{\sbl k}_{\sbl l}(\bla,\blc,m)
\Big(\left\langle \sR_n, B^m_{\sbl l} \right\rangle_{\sbla}
-\left\langle \sT_m, B^m_{\sbl l} \right\rangle_{\sbla}\Big),\qquad \bl k\in\Omega^{\sbl c}_m,
\]
where we have used Lemma~\ref{L:bestpol}.
We obtain
\begin{align*}
	\left\langle \sR_n, \,B^m_{\sbl l}\right\rangle_{\sbla}
	  =& 
	           \sum_{\sbl h\in\Theta_n}\omega_{\sbl h}r_{\sbl h}
	           \left\langle \frac{B^n_{\sbl h}}{\omega},\,B^m_{\sbl l}\right\rangle_{\sbla}\\
	  =& 
	           \sum_{\sbl h\in\Theta_n}\omega_{\sbl h}r_{\sbl h}
	           \binom{n}{\bl h}\binom{m}{\bl l}\rbinom{n+m}{\bl h+\bl l}
	           \left\langle \frac{1}{\omega},\,B^{n+m}_{\sbl h+\sbl l}\right\rangle_{\sbla}\\
	  =& 
	           \sum_{\sbl h\in\Theta_n}\omega_{\sbl h}r_{\sbl h}
	           \binom{n}{\bl h}\binom{m}{\bl l}\rbinom{n+m}{\bl h+\bl l}\,I_{\sbl h+\sbl l},	                    
\end{align*}
where we use the notation \eqref{E:I}. 
Further, using  equations \eqref{E:Ber2} and \eqref{E:Jinprod}, we obtain
\begin{align*}
	\left\langle \sT_m, B^{m}_{\sbl l} \right\rangle_{\sbla}=  &
	\sum_{\sbl h\in\Gamma^{\sbl c}_m}g_{\sbl h}
	        \left\langle B^m_{\sbl h}, B^{m}_{\sbl l} \right\rangle_{\sbla}\\	
	        = &\sum_{\sbl h\in\Gamma^{\sbl c}_m}g_{\sbl h}	     
	       \binom{m}{\bl h}\binom{m}{\bl l}
	       \frac{(\alpha_1+1)_{h_1+l_1}(\alpha_2+1)_{h_2+l_2}(\alpha_3+1)_{2m-|\sbl h|-|\sbl l|}} {(|\bla|+3)_{2m}}. 
\end{align*}
Hence, the formula \eqref{E:psc} follows.
\end{pf}

\begin{rem}\label{R:integr}
	In general, the integrals \eqref{E:I} cannot be evaluated exactly. In Section~\ref{SS:impl-I}, 
	we  show that they can be efficiently computed numerically up to high precision 
	using an extension of the method of \cite{Kel07}.
	
	  In the special case where all the weights  $\omega_{\sbl i}$, $\bl i\in\Theta_n$, are equal,
	 the rational function
	 \eqref{E:Rsc} reduces to a polynomial of degree $n$, so that  
	the  problem is actually the constrained polynomial degree reduction problem (see, e.g., \cite{WL10}).
	Evaluation of the integrals is then a~simple task.
	
\end{rem}

\section{Implementation of the method}
						\label{S:impl}
In this section, we discuss some computational  details  of  the polynomial approximation of the rational
B\'ezier function, described in Section~\ref{S:main} (see Theorem~\ref{T:main}).

\subsection{Computing the coefficients $E^{\sbl k}_{\sbl l}(\bla,\blc,m)$}
						\label{SS:impl-E}

We have to compute all the coefficients $E^{\sbl k}_{\sbl l}(\bla,\blc,m)$ with ${\bl k},\,\bll\in\Omega^{\sbl c}_m$. It has been shown \cite{LKW15} that they can be given   in terms of 
\[ 		
	e^{\sbl k}_{\sbl l}(\bl\mu,M):=\langle D^M_{\sblk},D^M_{\sbll}\rangle_{\sbl\mu},
	\qquad\blk,\bl l\in\Theta_M,
\] 	
with $M:=m-|\blc|$ and $\blm:=\bla+2\blc$, where  $D^M_{\sblk}\equiv D^M_{\sbl k}(\bl\cdot;\bl\mu)$ 
are the unconstrained dual Bernstein polynomial of total degree $M$ (cf. \eqref{E:dBer2}). 
See Eq. \eqref{E:E-e} for details.
Obviously, 
\(
	e^{\sbl k}_{\sbl l}(\bl\mu,M)=e^{\sbll}_{\sblk}(\bl\mu,M).
	\)
The following algoritm  is  based on the  recurrence relations satisfied by
 $e^{\sblk}_{\sbll}\equiv e^{\sbl k}_{\sbl l}(\bl\mu,M)$, obtained in  the paper cited above.
\begin{alg}[Computing the coefficients $E^{\sbl k}_{\sbl l}
(\bla,\blc,m)$]
	\label{A:Ecomp}%
	\   
\begin{description}
\itemsep4pt
\item[{\sc Step 1}] Let $M:=m-|\blc|$, $\blm:=\bla+2\blc$.	
\item[{\sc Step 2}] For $l_1=0,1,\ldots, M-1$,\\
\hspace*{4.0em}$l_2=0,1,\ldots,M-l_1$,\\   
\hspace*{2.2em}compute  
\begin{equation}\label{E:E00}
	e^{\sbl 0}_{\sbl l} :=\frac{(-1)^{l_1}(|\blm|+3)_
{M}}{M!(\mu_1+2)_{l_1}}
	\sum_{i=0}^{M-l_1}C^\ast_i\,h_i(l_2;\mu_2,\mu_3,M-
l_1),
	\end{equation}
\hspace*{2.2em}where $\bl 0=(0,0)$, $\bl l=(l_1,l_2)$,  
$h_i(t;a,b,N)$ are the Hahn polynomials (cf.~\eqref
{E:Hahn1}), and
\[
C^\ast_i:=\left\{\begin{array}{ll}
\dfrac{(\mu_1+2)_{M}}{(|\blm|-\mu_1+2)_{M-l_1}},&\quad i=0,
\\[2.5ex]
		\displaystyle (-1)^{i}
		\dfrac{(2i+|\blm|-\mu_1+1)(\mu_1+2)_{M-i}(|
\blm|+M+3)_{i}}
		{i!(\mu_3+1)_i(|\blm|-\mu_1+i+1)_{M-
l_1+1}},&\quad i\ge1;	
			\end{array}\right.
\]
\hspace*{2.2em}next put $e^{\sbl l}_{\sbl 0}:=e^{\sbl 0}_
{\sbl l}$.
   
\item[{\sc Step 3}] For $k_1=0,1,\ldots,M-1$,
\begin{description}
\itemsep4pt
\item[$1^o$]
	for $k_2=0,1,\ldots,M-k_1-1$,\\
	\hspace*{1.55em}$l_1=k_1,k_1+1,\ldots,M$,  \\
	\hspace*{1.55em}$l_2=0,1,\ldots,M-l_1$,  \\
    compute 
    \[ 
	e^{\sbl k+\sbl v_2}_{\sbl l}:=\left([\sigma_1(\bl 
k)-\sigma_1(\bl l)]e^{\sbl k}_{\sbl l}
			-\sigma_2(\bl k) e^{\sbl k-\sbl 
v_2}_{\sbl l} 
			+\sigma_0(\bl l) e^{\sbl k}_{\sbl l
+\sbl v_2}
			+\sigma_2(\bl l) e^{\sbl k}_{\sbl 
l-\sbl v_2}\right)/\sigma_0(\bl k),	
    \] 
    where $\bl k=(k_1,k_2)$, $\bl l=(l_1,l_2)$, $\bl v_2:=
(0,1)$, and where for $\bl t:=(t_1,t_2)$ we define
\[ 
	\;\;\sigma_0(\bl t):=(|\bl t|-M)(t_2+\mu_2+1),\quad
\!
	\sigma_2(\bl t):=t_2(|\bl t|-\mu_3-M-1), \quad\!
	\sigma_1(\bl t):=\sigma_0(\bl t)+\sigma_2(\bl t),
\] 
next put $e^{\sbl l}_{\sbl k+\sbl v_2}:=e^{\sbl k+\sbl v_2}
_{\sbl l}$;
\item[$2^o$]for $l_1=k_1+1,k_1+2,\ldots,M$,  \\
	\hspace*{1.55em}$l_2=0,1,\ldots,M-l_1$,  \\
    compute 
    \[
	e^{\sbl k+\sbl v_1}_{\sbl l}
	:=\left([\tau_1(\bl k)- \tau_1(\bl l)]e^{\sbl k}_
{\sbl l}
				  -\tau_2(\bl k)\,e^{\sbl 
k-\sbl v_1}_{\sbl l}
				 + \tau_0(\bl l)\,e^{\sbl 
k}_{\sbl l+\sbl v_1}				  
				 + \tau_2(\bl l)\,e^{\sbl 
k}_{\sbl l-\sbl v_1}\right)/\tau_0(\bl k),		
    \] 
where  $\bl k=(k_1,0)$, $\bl l=(l_1,l_2)$, $\bl v_1:=(1,0)$, 
 and for   $\bl t:=(t_1,t_2)$ the coefficients  
$\tau_{j}(\bl t)$ are given by
\[ 
	\tau_0(\bl t):=(|\bl t|-M)(t_1+\mu_1+1), \quad\!	
	\tau_2(\bl t):=t_1(|\bl t|-\mu_3-M-1),  \quad\!
	\tau_1(\bl t):=\tau_0(\bl t)+\tau_2(\bl t);
\] 
next put $e^{\sbl l}_{\sbl k+\sbl v_1}:=e^{\sbl k+\sbl v_1}
_{\sbl l}$.     
\end{description} 
\item[{\sc Step 4}] For $ \blk,\,\bll\in\Omega^{\sbl c}_m$, 
compute 
	\begin{equation}
		\label{E:E-e}
			E^{\sbl k}_{\sbl l}(\bla, c,m):=U
\,V_{\sbl k}\,V_{\sbl l}\,
	\,e^{\sbl k-\sbl c'}_{\sbl l-\sbl c'},
	\end{equation}
	\hspace*{2.2em}where $\bl c':=(c_1,c_2)$, and
\[
	 U:=
	(|\bla|+3)_{2|\sbl c|}\prod_{i=1}^{3}
	(\alpha_i+1)_{2c_i}^{-1},\qquad
	V_{\sbl h}:= \binom{m-|\bl c|}{\bl h-\bl c'}\rbinom
{m}{\bl h}.
		\]
\end{description}
\end{alg}

As noticed in Remark~\ref{R:Clenshaw},  the sum
in \eqref{E:E00}  can be evaluated  efficiently  using the Clenshaw's algorithm, 
at the cost of $O(M-l_1)$ operations.

\subsection{Computing the integrals $I_{\sbl j}$}
										\label{SS:impl-I}

The most computationally expensive part of the proposed method is
the evaluation of the collection of integrals \eqref{E:I}. For example,
for $n+m=22$, if $\bl c=(0,0,0)$, there are $276$ two-dimensional integrals
to be computed. It is obvious that using any standard quadrature would
completely ruin the efficiency of the algorithm. Moreover, if any of the
parameters $\alpha_i$ ($i=1,2,3$) in (\ref{E:w}) is smaller
than $0$ and the corresponding constrain parameter $c_i$
equals zero, then the integrands in (\ref{E:I}) are singular
functions, and standard quadratures may fail to
deliver any approximations to the integrals.

Therefore, for evaluating the complete set of integrals (\ref{E:I}),
we introduce a special scheme which is based on the general method \cite{Kel07}
for approximating singular integrals. The proposed numerical quadrature is of
the automatic type, which means that the required number of nodes is
adaptively selected, depending on the complexity of the rational
B\'ezier function, so that the requested accuracy of the
approximation is always achieved. Most importantly,
the algorithm is extremely effective in the considered application.
In the example given at the beginning of this subsection ($n+m=22$), the
time required to compute the whole collection of $276$ integrals is only
twice\footnote{Based on the Maple implementation of the algorithm.
If the collection consists of 990 integrals ($n+m=42$), the computation
time increases by only 50\% (compared to the case of 276 integrals).
The detailed report from the efficiency test can be found at the
end of Appendix B
.} longer than the time needed to approximate a single
separate integral of a similar type.

First, we shall write the integral (\ref{E:I}) in a different form
which is better suited for fast numerical evaluation. Observe that bivariate
Bernstein polynomials \eqref{E:Ber2} can be expressed in terms of
univariate Bernstein polynomials. Namely, we have
\[
  B^N_{\sbl j}(\bl x)= B^N_{j_1}(x_1) B^{N-j_1}_{j_2}({x_2}/{(1-x_1)}),
  \qquad \bl j:=(j_1,j_2),\; \bl x:=(x_1,x_2),
\]
where $B^M_i(t):=\binom{M}{i}t^i(1-t)^{M-i}$,  $0\le i\le M$, are univariate Bernstein polynomials.
Further, the bivariate weight function $w_{\sbla}$  (see \eqref{E:w}) can be expressed as
\[
  w_{\sbla}(\bl x)
  = A_{\sbla}\,v_{\alpha_2+\alpha_3,\alpha_1}(x_1)
    \,v_{\alpha_3,\alpha_2}({x_2}/{(1-x_1)}),
\]
where $v_{\alpha,\beta}(t):=(1-t)^\alpha t^\beta$ is the univariate Jacobi weight function.
Hence, the integral \eqref{E:I} can be written as
\begin{align}
\nonumber
  I_{\sbl j}&=\int_{0}^{1}\int_{0}^{1-x_1} w_{\sbla}(\bl x)\frac{B^{N}_{\sbl j}(\bl x)}{\omega(\bl x)}\,{\rm d}x_2\,{\rm d}x_1 \\
\nonumber
  &=A_{\sbla}\int_{0}^{1}v_{\alpha_2+\alpha_3+1,\alpha_1}(s)B^{N}_{j_1}(s)
  \left(\int_{0}^{1}v_{\alpha_3,\alpha_2}(t)\frac{B^{N-j_1}_{j_2}(t)}{\omega^\ast(s,t)}\,{\rm d}t\right){\rm d}s \\[0.5ex]
\label{E:I2}
  &=A_{\sbla}\binom{N}{\bl j} 
  \int_{0}^{1}v_{a,b}(t)
  \left(\int_{0}^{1}v_{c,d}(s)\frac{1}{\omega^\ast(s,t)}\,{\rm d}s\right){\rm d}t,
\end{align}
where we denoted $N:=n+m$, 
\begin{equation}\label{E:abcd}
\left.\begin{array}{ll}
  a\equiv a(\bl j):=\alpha_3+N-|\bl j|,&\quad  b\equiv b(j_2):=\alpha_2+j_2,\\[1ex]
  c\equiv c(j_1):=\alpha_2+\alpha_3+N-j_1+1, &\quad
  d\equiv d(j_1):=\alpha_1+j_1,
\end{array}\,\right\}
\end{equation}
and
\begin{equation}
\label{E:omegastar}
  \omega^\ast(s,t):=\omega(s,(1-s)t)
  =\sum_{i=0}^n w_i(t) B^n_{i}(s), \qquad
  w_i(t) = \sum_{j=0}^{n-i}\omega_{i,j}\,B^{n-i}_{j}(t).
\end{equation}
Note that the computation of values of the integrand is now much
more effective, because the coefficients $w_i$ of the function $\omega^\ast$
($1\leq i\leq n$) in (\ref{E:omegastar}) do not depend of the inner integration
variable $s$. The main idea is, however, to compute the values of $\omega^\ast$
only once (at a properly selected set of quadrature nodes), and obtain a tool
for fast computation of the integrals (\ref{E:I2}) for different values
of $a$, $b$, $c$, and $d$, i.e.\ for different values of $\bl j$.

For arbitrary fixed $t\in[0,1]$, define the function
\begin{equation}\label{E:psi}
  \psi_t(s):=\omega^\ast(s,t)^{-1}.
\end{equation}
It is easy to see that we can write
\[
  I_{\sbl j}=A_{\sbla}\binom{N}{\bl j} J(a,b,\Phi),
\]
with
\begin{equation}\label{E:Phi}
  \Phi(t):=J(c,d,\psi_t),
\end{equation}
  where we use the notation
\[
  J(\alpha,\beta,f):=\int_{0}^{1}(1-x)^\alpha x^\beta f(x){\rm d}x.
\]
The functions $\psi_t$ and $\Phi$ are analytic in a closed complex
region containing the interval $[0,1]$ (it is proved in Appendix B). 
This implies that (cf.\ \cite[Chapter 3]{Riv90}) they can be accurately and
efficiently approximated by polynomials given in terms of
the (shifted) Chebyshev polynomials of the first kind,
\begin{equation}\label{E:TchExpan}
\begin{array}{l}
\displaystyle \psi_t(s)\approx S_{M_t}(s) := \sum_{i=0}^{M_t}{'\,}\gamma^{[t]}_{i} T_i(2s-1),\\[3ex]
\displaystyle \Phi(t) \approx \hat{S}_{M}(t) :=
    \sum_{l=0}^{M}{'\,}\hat{\gamma}_l T_i(2t-1),
    \end{array}\qquad 0\leq s,t\leq 1,
\end{equation}
where $M$ may depend on $j_1$, and the prime denotes a sum with the first term halved.
Once the above expansions are computed (this can be done in a time proportional
to $M_{t}\log(M_{t})$ and $M\log(M)$), the integrals $J(\cdot,\cdot,\cdot)$
can be easily evaluated using the following algorithm
that was proved in \cite{LWK12}.

\begin{alg}[Computing the integral $ J(\alpha,\beta;S)$, $S$ being a polynomial] \label{A:PK_1d}
\  \\
Given numbers $\alpha,\,\beta>-1$, let $r := \beta-\alpha$, $u := \alpha+\beta+1$.
Let  $S_{\MM}$ be a polynomial defined by
\[
  S_{\MM}(x) = \sum\limits_{i=0}^{\MM}{'} \gamma_i T_i(2x-1).
\]
Compute the sequence $d_i$, $0\leq i \leq {\MM}+1$,  by
\begin{align*}
  &d_{{\MM}+1} = d_{\MM} := 0, \\[1.0ex]
  &d_{i-1} :=
    \frac{2 r d_{i} + (i-u) d_{i+1} - \gamma_{i}}
         {i+u},
    \qquad i={\MM},{\MM}-1,\dots,1.
\end{align*}
\noindent \textbf{Output}:
$ J(\alpha,\beta;S_{\MM}) = \BB\cdot \left( \frac12\gamma_0 - r d_0 + u d_1 \right)$,
where $\hts\BB := \Gamma(\alpha+1)\Gamma(\beta+1) / \Gamma(\alpha+\beta+2)$.
\end{alg}

By the repeated use of the above very fast scheme, we may efficiently approximate
the whole set of integrals $I_{\sbl j}$ for $\bl j\in\Omega_{n+m}^{\sbl c}$.
The remaining technical details of the adaptive implementation of
the proposed quadrature and the complete formulation of the
integration algorithm are presented in Appendix A.

\subsection{ Main algorithm}
						\label{SSS:MainAlg} 

The method  presented in this paper is summarized in the following algorithm.
\begin{alg}[Polynomial approximation of the rational B\'ezier triangular surface]	
	\label{A:RTBSappr}
Given the coefficients $r_{\sbl k}$ 
and positive weights $\omega_{\sbl k}$, $\bl k\in\Theta_n$, of the rational 
function \eqref{E:Rsc},
the coefficients $p_{\sbl k}$ of the degree $m$ polynomial \eqref{E:Psc}, minimising the error
\eqref{E:er-sc},
with the constraints \eqref{E:gsc},
can be computed in the following way.\\[-3ex]
\begin{description}
\itemsep0.25pt
\item[{\sc Step 1}] Compute the  table 
	 $\{E^{\sbl k}_{\sbl l}(\bla, c,m)\}_{\sbl k,\sbl l\in\Omega_m^{\sbl c}}$ 
	by Algorithm~\ref{A:Ecomp}. 
\item[{\sc Step 2}] Compute the table $\{I_{\sbl j}\}_{\sbl j\in\Omega^{\sbl c}_{n+m}}$ by  Algorithm~\ref{A:I-comp}.
   
\item[{\sc Step 3}] For $\bl k\in\Gamma^{\sblc}_m$, put $p_{\sbl k}:=g_{\sbl k}$.
  
\item[{\sc Step 4}] For $\bl k\in\Omega^{\sblc}_m$, compute $p_{\sbl k}$ by \eqref{E:psc}.
\end{description}
\noindent \textbf{Output}: Set of the coefficients $p_{\sbl k}$, $\bl k\in\Theta_m$.
\end{alg}


\section{Examples}
					\label{S:exmp}

In this section, we present some examples of approximation of rational  triangular B\'ezier patches
by triangular B\'ezier patches. 
No theoretical justification is known for the best choice of the vector parameter
$\bla$ in the distance functional \eqref{E:dist} if we use the \textit{error function} 
\begin{equation}
	\label{E:err}
	\Delta(\bm x):=\|\sR_n(\bm x)-\sP_m(\bm x)\|		
	\end{equation}
to measure the quality of the approximation. On the base of numerical experiments, we claim that
$\bla=(-\frac12,-\frac12,-\frac12)$  usually leads to slightly better results than the ones obtained for
other "natural" choices of parameters, including the usually preferred $\bla=(0,0,0)$
(meaning $w_{\sbla}(\bl x)=1$).  The computations were performed in 
16-decimal-digit  arithmetic. In the implementation of  Algorithm~\ref{A:I-comp},   we have assumed $\varepsilon=5\times10^{-16}$ in \eqref{E:stop}, and used the initial values $M^{*}=M_{k}^{*}=32$. 

\subsection{Example 1}\label{SS:example1}

Let $\sR_6$ be the degree 6 rational triangular B\'ezier patch
\cite[Example 2]{Hu13},
\begin{equation}
	\label{E:R6}
	\sR_6(\bl x):=
	\dfrac{\displaystyle \sum_{\sbl k\in\Theta_6}\omega_{\sbl k}r_{\sbl k}B^6_{\sbl k}(\bl x)}
	{\displaystyle \sum_{\sbl k\in\Theta_6}\omega_{\sbl k}B^6_{\sbl k}(\bl x)}, \qquad \bl x\in T,	
\end{equation}
$T$ being the standard triangle \eqref{E:T}, and the control points $r_{\sbl k}$ and the associated
weights $\omega_{\sbl k}$ being listed in Table~\ref{tab:R6cp+w}. 
\begin{table}[phtb]
	\caption{Control points $r_{\sbl k}$ (\textit{upper entries}) and the associated
		weights $\omega_{\sbl k}$ (\textit{lower entries}) of the surface \eqref{E:R6},
	with $\bl k=(k_1,k_2)\in\Theta_6$.\label{tab:R6cp+w}}
{\footnotesize
\[
	\begin{array}{c|ccccccc}
	k_1\setminus k_2  &  0 & 1 & 2 & 3 & 4 &5 &6 \\[1ex]\hline 
	&&&&&&&\\
0 & (6,0,2) & (5,0,3) & (4,-0.5,3.5) &
                       (3,-0.2,4) & (1.5,0.5,2) & (0.4,0.4,1)&(0,0,0)\\  
   &0.8 & 0.3 & 1.8 & 1.2 &0.8 & 0.2&1.6\\                         
1 & (5.2,1,3)& (4.5,1,3) & (3,0.6,4) &
                                        (2,0.9,3)  &  (1.2,1,2)  & (0.4,0.8,0.6)&\\ 
   &1& 0.4 &0.8&2.4 &1.3 &0.9& \\                                     
2 & (4.5,2,5) & (4,2.2,4) & (3,2,3)    & (2,1.2,2)    &(0.8,1.5,1.5)&& \\ 
   & 0.5& 1 & 1  & 1.8 &0.8&& \\ 
3 & (4,3,6) & (2.5,2.5,5) & (1.5,2.8,4)   &(1,2,3)&& &\\ 
   & 0.3 & 2 & 1 &0.9&&& \\ 
4 & (3.5,3.5,4) & (2.5,3,5) &(1.5,3.5,3)&&&&\\ 
   & 1.5& 0.6 &1.2&&& &\\ 
5 &  (3,4.2,2)&(2,4,2)  &  &  &  & &\\
   & 0.8 & 0.5 &  &  &  & &\\
6 &  (2,5,1)&  &  &  &  & &\\
   & 1 &  &  &  &  & &   
\end{array}
\]
}
\end{table}
We let $\bla=(-\frac12,-\frac12,-\frac12)$, $\blc=(1,1,1)$, and constructed the degree 5 
best approximating polynomial patch
	\[
		\sP_5(\bl x):=\sum_{\sbl k\in\Theta_5}p_{\sbl k}B^5_{\sbl k}(\bl x),
		 \qquad \bl x\in T,
	\]
	under the restriction
	 	$p_{\sbl k}=g_{\sbl k}$  for $\bl k\in\Gamma^{\sbl c}_5$, 	
	 	where
	 	\[
	 		\Gamma^{\sbl c}_5:=\{\bl k=(k_1,k_2):k_1=0,\; \mbox{or}\; k_2=0,
	 	 \; \mbox{or}\; |\bl k|=5\},
	 	\]
         and the set of points $g_{\sbl k}$, $\bl k\in\Gamma^{\sbl c}_5$, is obtained 
         in the following way. As well known, the boundary of the patch \eqref{E:R6} 
         is formed by  three degree 6  rational  B\'ezier curves.      
         The  least squares degree 5 polynomial approximation to each of these rational curves,
         with the endpoints preservation, is constructed
         using an extension of  the method of \cite{LWK12}, described in \cite{LWK15} 
         (the input data: $m=5$, $\alpha=\beta=-\frac12$, $k=l=1$, 
         	 notation used being that of \cite{LWK15}). 
         Now, the set of points $g_{\sbl k}$ is the appropriate collection of all control
         points of the three resulting B\'ezier curves.
         
\begin{figure}[htb]
	\begin{center}  
		
		\includegraphics[scale=0.32,angle=-90]{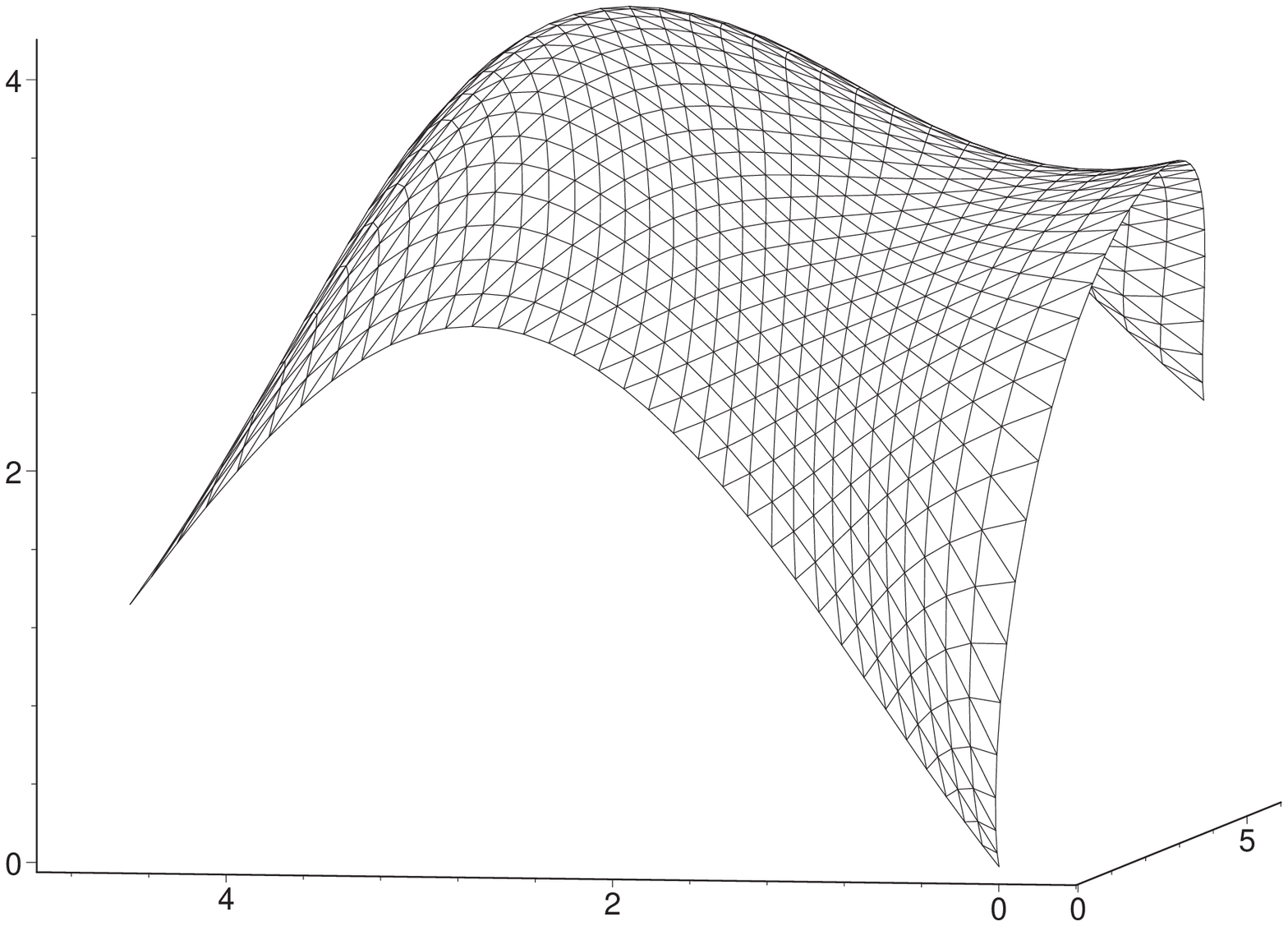}%
	\hspace*{0.8cm}%
	{\includegraphics[scale=0.32,angle=-90]{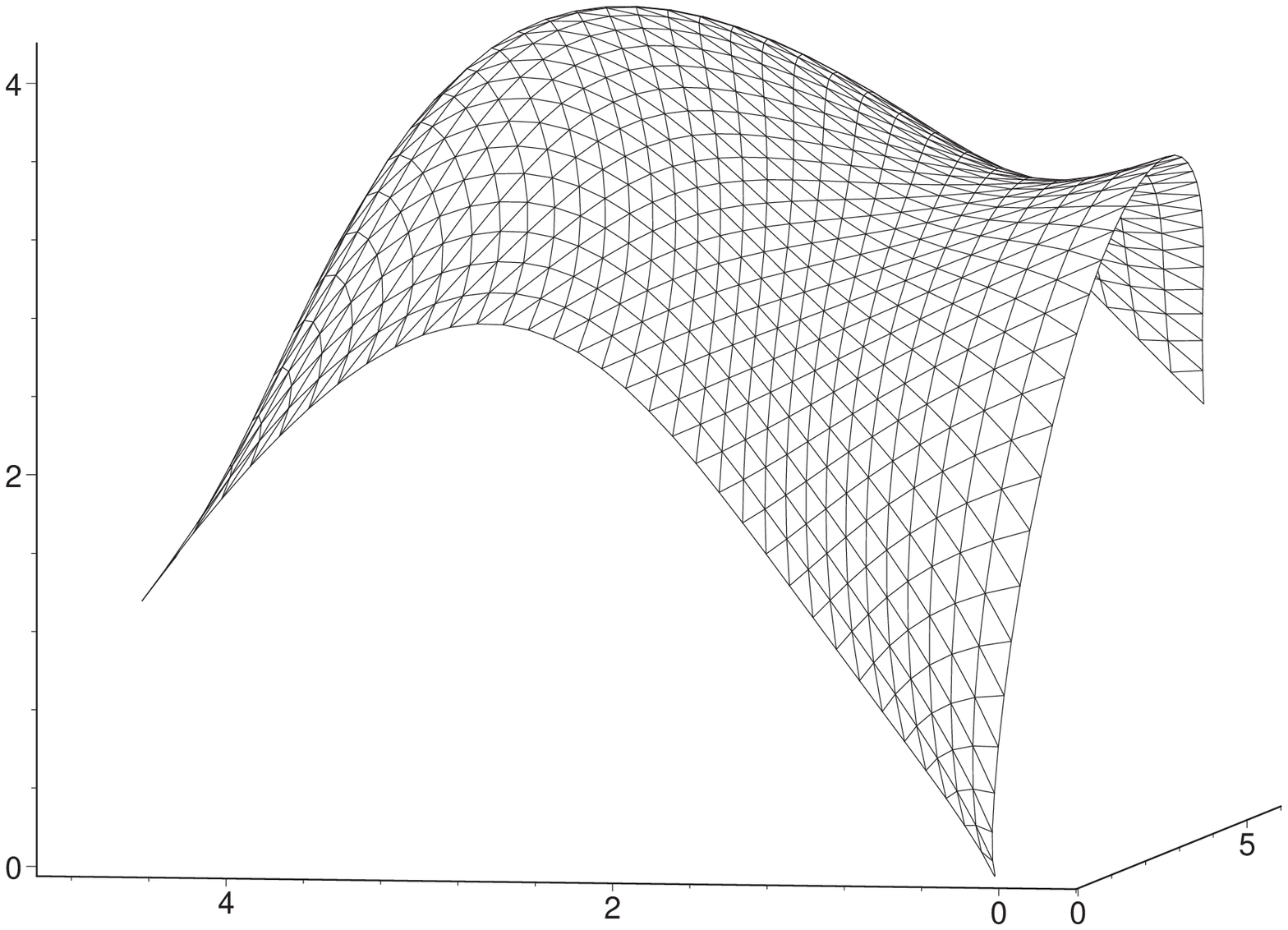}}%
	
	\vspace*{0.5cm}%
	\hspace*{-0.32cm}%
	\includegraphics[scale=0.3,angle=-90]{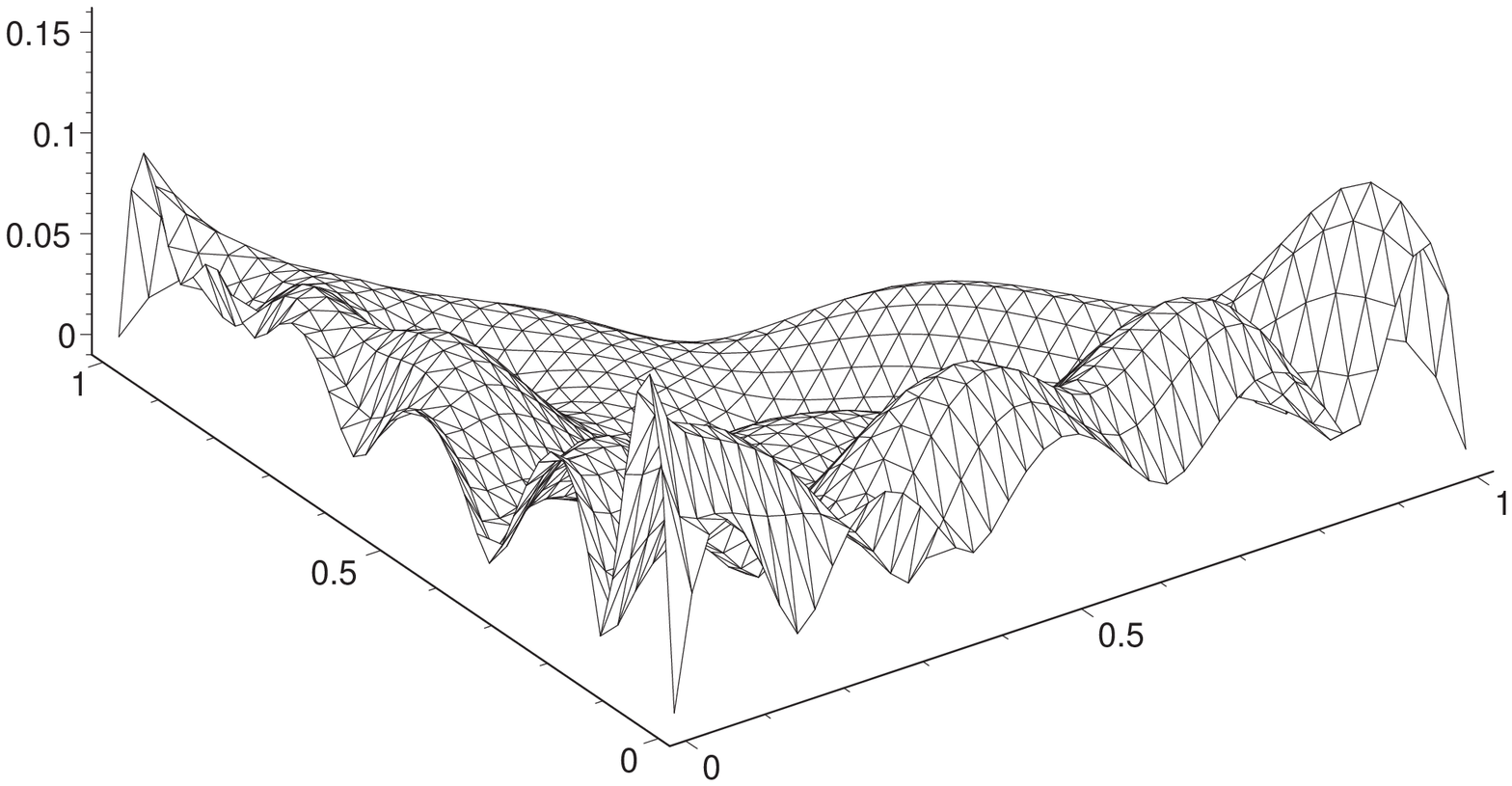}%
	\hspace*{0.8cm}%
	{\includegraphics[scale=0.3,angle=-90]{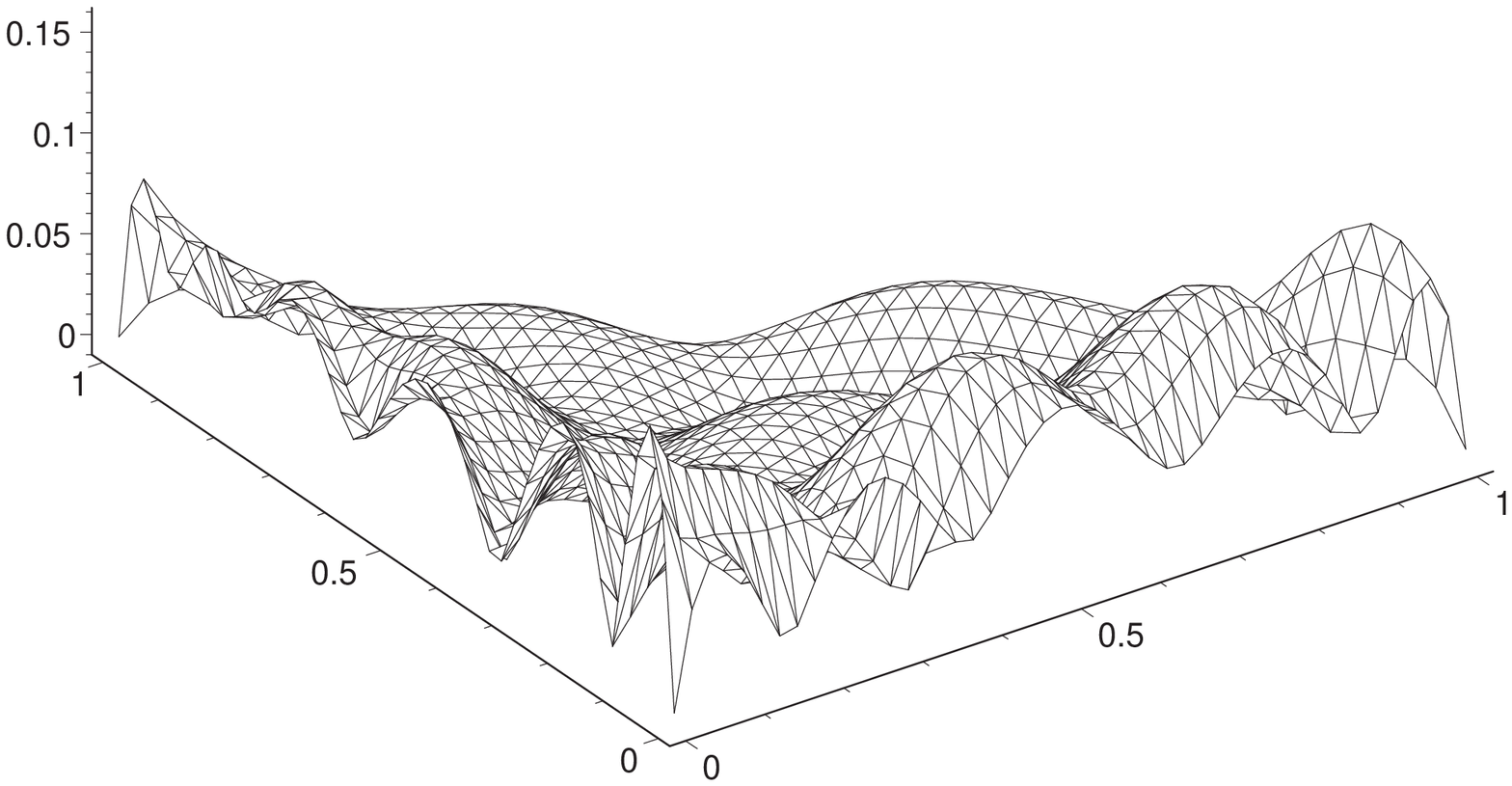}}%
	\caption{\small Constrained degree 5 polynomial approximation of the degree 6 rational triangular 
		B\'ezier surface, with  $\blc=(1,1,1)$. \textit{Upper part}:  Rational surface $R_6(\bm x)$ and the approximating surface $\sP_5(\bm x)$
		with $\bla=(-\frac12,-\frac12,-\frac12)$. \textit{Lower part}: The error $\Delta(\bm x)$ plots corresponding to $\bla=(0,0,0)$  and $\bla=(-\frac12,-\frac12,-\frac12)$, respectively.
		The maximum errors  are $0.16$ and $0.13$, respectively.		 	          
	 Notice that 
	 the original surface and the approximating surface 
	 agree at the corner points.%
	\label{fig:Fig2}}
	\end{center}
\end{figure}

         We have repeated the computations for  $\bla=(0,0,0)$ 
         (with $\alpha=\beta=0$, in \cite{LWK15}), obtaining slightly worse results.
         The maximum errors $\max_{\sbl x\in T}\Delta(\bm x)$ (cf. \eqref{E:err}) of the obtained results (see Fig.~\ref{fig:Fig2}) are about 50\% less than those reported in \cite[Table 1]{Hu13}.

\subsection{Example 2}\label{SS:example2}
\begin{figure}[phtb]
	\begin{center}
        \includegraphics[scale=0.52,angle=-90]{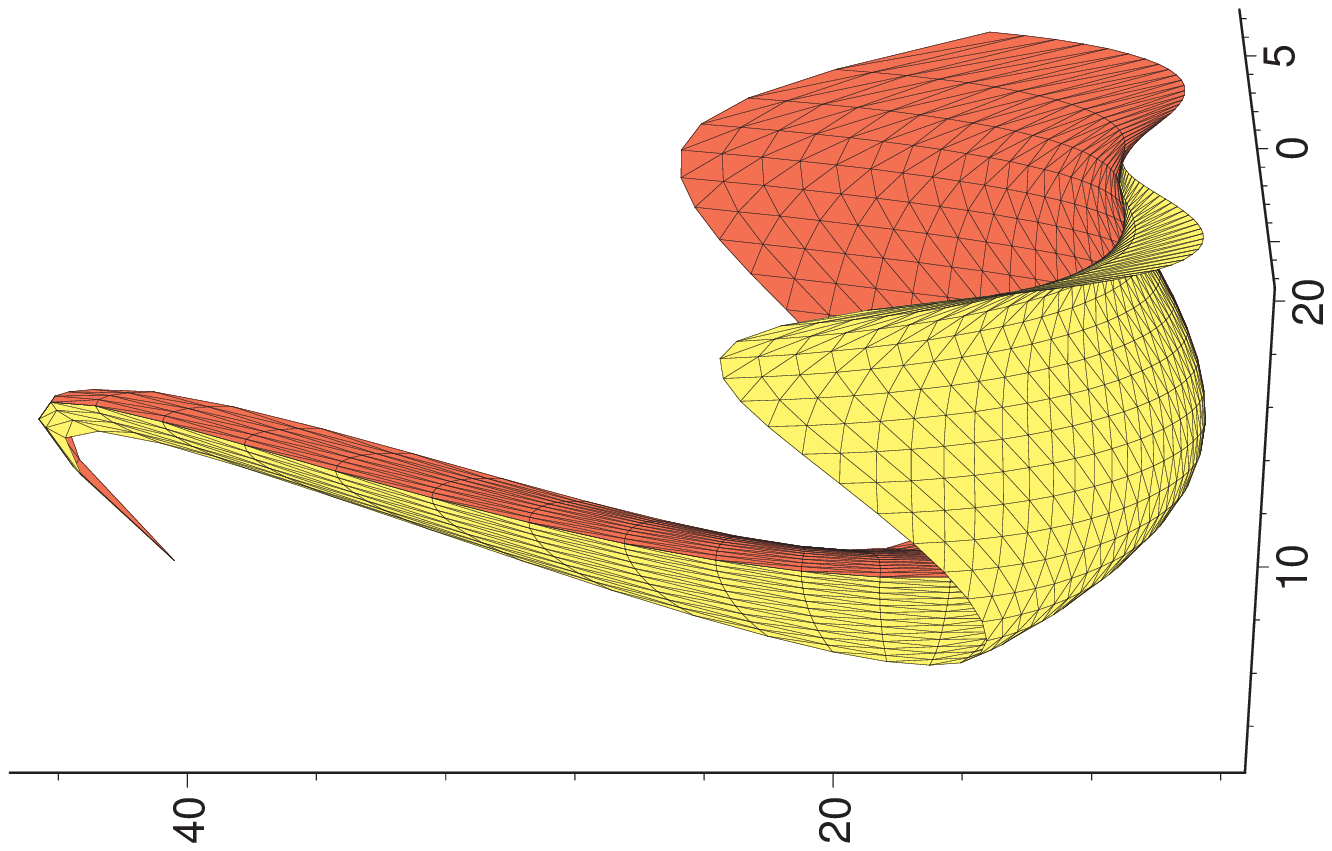}%
         \hspace*{0.3cm}
         \includegraphics[scale=0.52,angle=-90]{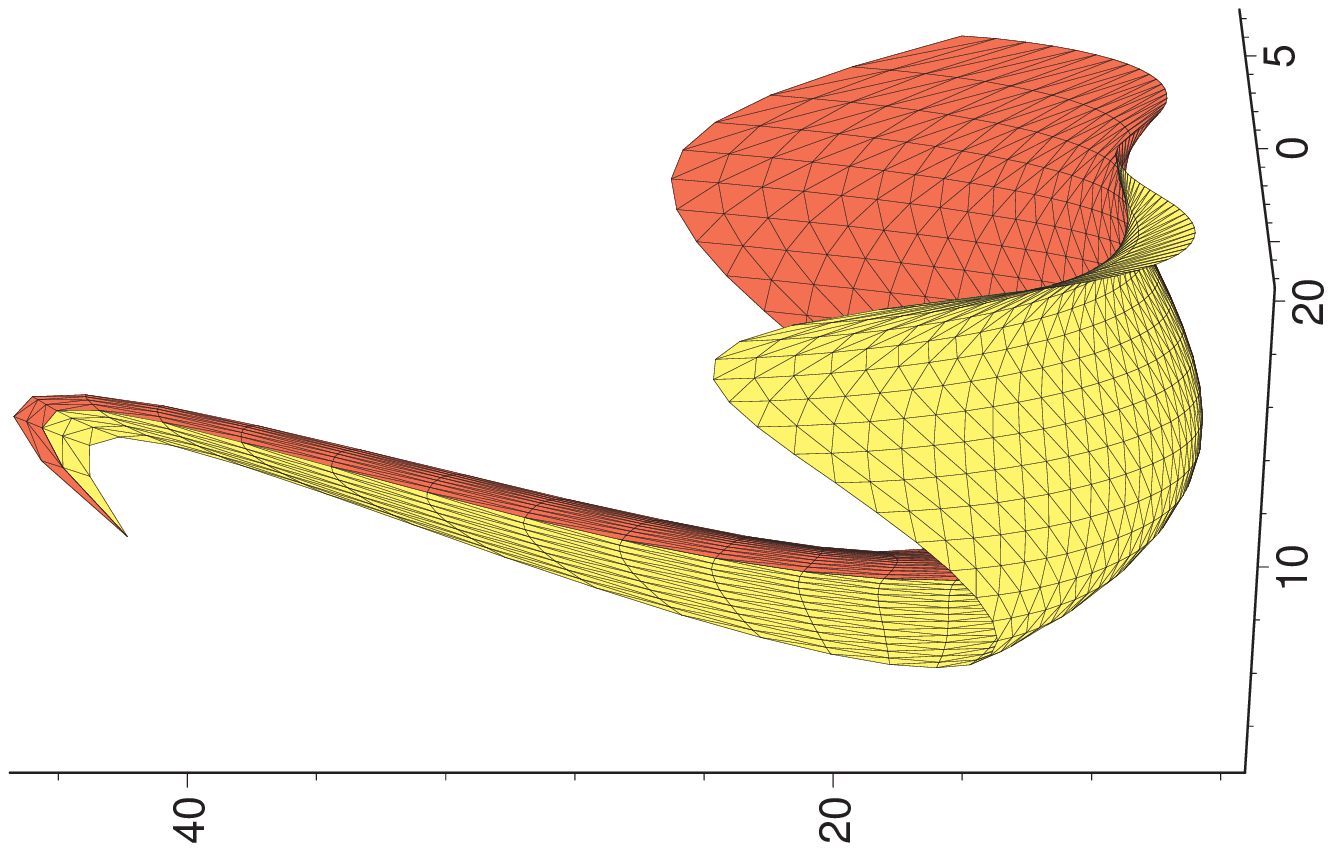}%
        \hspace*{0.3cm}
	\includegraphics[scale=0.52,angle=-90]{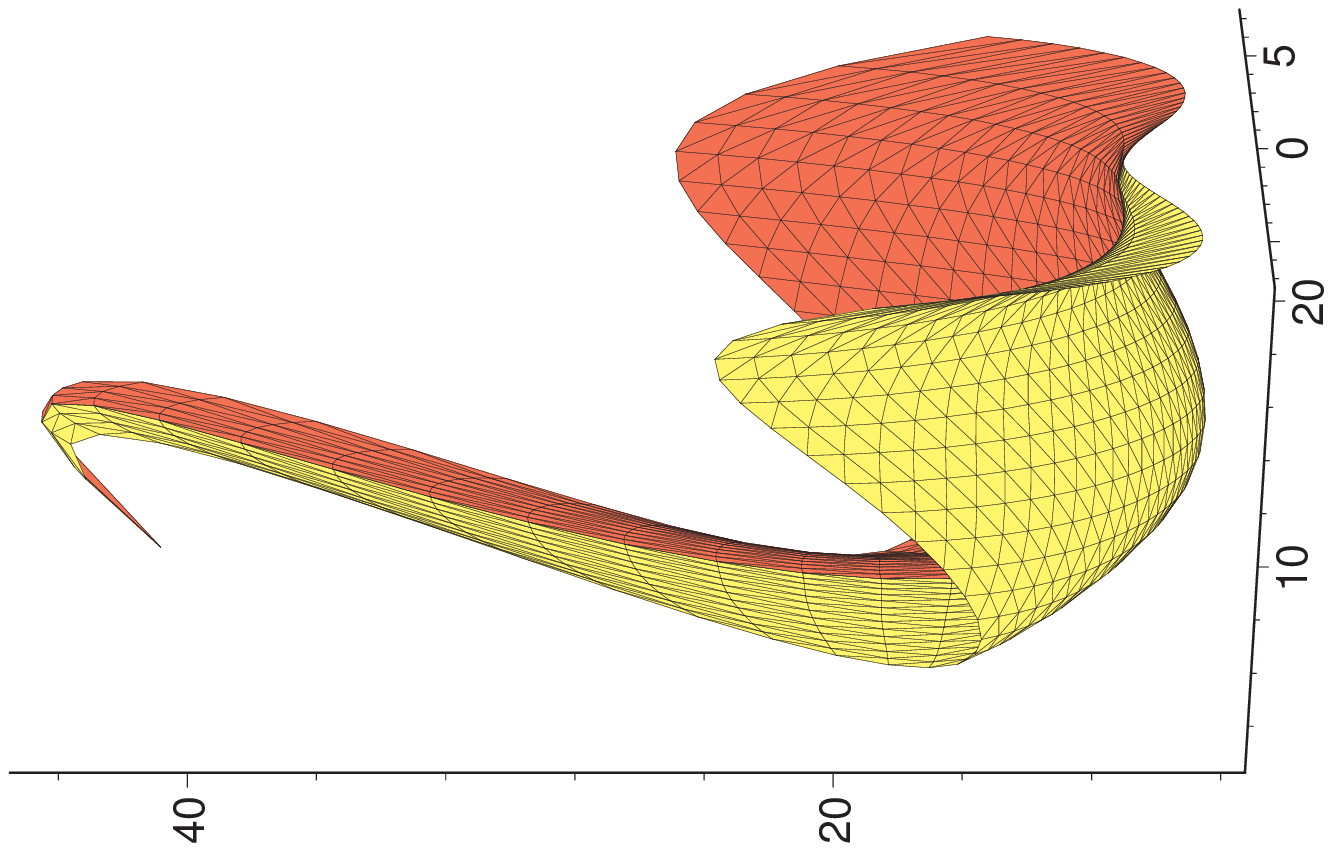}%
	\caption{\small  The composite rational B\'ezier surface \eqref{E:Rast}
		(\textit{left}) and the  $C^1$-continuous composite polynomial surfaces of degree (5,5) (\textit{middle}) and (6,6) (\textit{right}).	          
	\label{fig:Fig3}}
	\end{center}
	\end{figure}

Let $\sR^\ast$ be the composite  rational  surface,
\begin{equation}
	\label{E:Rast}
	\sR^\ast(\bl x):=\left\{\begin{array}{ll}
	\sR^R_5(\bl y), &\qquad\bl y:=(1-|\bl x|,x_1-x_2),\;\bl x\in T_R,\\[2ex]
	\sR^Y_5(\bl z), &\qquad\bl z:=(x_2-x_1,1-|\bl x|),\;\bl x\in T_Y,
	\end{array}\right.	
\end{equation}
where for $C\in\{R,Y\}$,
\begin{equation}
	\label{E:RC}
	\sR^C_5(\bl w):=
	\dfrac{\displaystyle \sum_{\sbl k\in\Theta_5}\omega^C_{\sbl k}r^C_{\sbl k}B^5_{\sbl k}(\bl w)}
	{\displaystyle \sum_{\sbl k\in\Theta_5}\omega^C_{\sbl k}B^5_{\sbl k}(\bl w)}, 
	\qquad \bl w\in T,	
	\end{equation}
$T$ being the standard triangle \eqref{E:T}, and
\begin{align*}
		T_R:=&\{\bl x=(x_1,x_2): x_1\ge x_2\ge0,\;|\bl x|\leq1\},  \\[1ex]
			  T_Y:=&\{\bl x=(x_1,x_2): x_2\ge x_1\ge0,\;|\bl x|\leq1\}.
		\end{align*}
	The control points $r^C_{\sbl k}$ and the associated
	weights $\omega^C_{\sbl k}$ of the rational patches \eqref{E:RC} can be found at the webpage
	\texttt{http://www.ii.uni.wroc.pl/\~{}pwo/programs.html}. 
	The surface \eqref{E:Rast} is shown in  Fig.~\ref{fig:Fig3} (the left plot).
	
        Now, we show how to obtain the degree $m$  polynomial approximations of the rational subpatches,
        which form a $C^1$-continuous composite surface.  
        
        $1^o$ Let $\sP^Y_m$ be the triangular B\'ezier 
        patch of degree $m$ approximating the rational patch $\sR^Y_5$ without constraints, i.e.,  for $\blc=(0,0,0)$.   
        Let $p^Y_{\sbl k}$ be the control points of the patch $\sP^Y_m$. 
	
        $2^o$ We approximate the rational patch $\sR^R_5$ by  the triangular B\'ezier 
        patch  $\sP^R_m$  of degree $m$, with constraints of the type $\blc=(2,0,0)$, 
        where the points $g_{\sbl k}\in\Gamma^{\sbl c}_m$ are chosen so that the $C^1$-continuity
        is obtained (cf. \cite[Section 17]{Far02}):
        \[
        	\begin{array}{ll}
        	g_{(0,i)}:= p^Y_{(i,0)},&\qquad i=0,1,\ldots,m,  \\[2ex]
        	g_{(1,i)}:=p^Y_{(i+1,0)}+(p^Y_{(i+1,0)}-p^Y_{(i,1)}),&\qquad i=0,1,\ldots,m-1.
        	\end{array}
        \]
	The results, obtained for $m=5$ and $m=6$, with $\bla=(-\frac12,-\frac12,-\frac12)$, 
	are shown in Fig.~\ref{fig:Fig3}. It can be observed that approximation of the rational 
	composite surface \eqref{E:Rast} by two jointed polynomial patches of  degree $m=5$ 
	(the middle plot) resulted in some visible differences.  Increasing the degree of the approximating polynomials to $m=6$ (the right plot) already
	gave a very satisfactory result. 
		
\section*{Appendix A: The adaptive algorithm for computing the integrals $I_{\sbl j}$}
						\label{S:AppA}
\setcounter{equation}0
\setcounter{subsection}0
\setcounter{thm}0
\renewcommand{\thethm}{A.\arabic{thm}}
\renewcommand{\thesection}{\Arabic{section}}
\renewcommand{\thesubsection}{A.\arabic{subsection}}
\renewcommand{\theequation}{A.\arabic{equation}}


We start with proving that the functions $\psi_t$ (\ref{E:psi}), $t\in[0,1]$, 
and $\Phi$ (\ref{E:Phi}) are analytic in a closed complex region containing the interval
$[0,1]$. The assertion is clearly true in the case of $\psi_t(z)=\omega^{*}(z,t)^{-1}$,
as the bivariate polynomial $\omega^{*}$ has no roots in $[0,1]\times[0,1]$.
Similarly, for any $s\in[0,1]$, the function $z\mapsto\omega^{*}(s,z)^{-1}$ is
analytic in a rectangular region $[-\sigma,1+\sigma]\times[-\sigma,\sigma]$,
where $\sigma>0$ does not depend on $s$. Thus, if $s\in[0,1]$, then
\begin{equation*}
  \int_{C} \omega^{*}(s,z)^{-1} {\rm d}z = 0
\end{equation*}
for any closed contour $C\subset[-\sigma,1+\sigma]\times[-\sigma,\sigma]$.
Consequently, if $\alpha,\beta > -1$, then
\begin{equation*}
  \int_{C} \Bigg( \int_{0}^{1} (1-s)^{\alpha}s^{\beta} \omega^{*}(s,z)^{-1} {\rm d}s\Bigg) {\rm d}z
  = \int_{0}^{1} (1-s)^{\alpha}s^{\beta} \Bigg( \int_{C} \omega^{*}(s,z)^{-1} {\rm d}z\Bigg) {\rm d}s
  = 0 .
\end{equation*}
Therefore, by Morera's theorem (see, e.g., \cite[Chapter 2.3]{Ahl79}),
the function $\Phi(z) = J(\alpha,\beta,\psi_z)$ 
is also analytic in $[-\sigma,1+\sigma]\times[-\sigma,\sigma]$.

The polynomials $S_{M_t}$ and $\hat{S}_M$ in (\ref{E:TchExpan}), which approximates
the functions $\psi_t$ and $\Phi$, are determined to satisfy the interpolation conditions
\[
  \left.\begin{array}{l}
    S_{M_{k}}(s_j) = \omega^\ast(s_j,t_k)^{-1},\quad  0\leq j\leq M_{k},\\[2ex]
    \hat{S}_{M}(t_k) = J(c,d;S_{M_{k}}),
  \end{array}\right\}\qquad 0\leq k\leq M,
\]
where, for simplicity, we denote $M_{k}\equiv M_{t_k}$,
and the interpolation nodes are given by
\begin{equation}
\label{E:sjtk}
  s_j = \frac{1}{2} + \frac{1}{2}\cos\frac{j\pi}{M_{k}},
    \qquad
  t_k = \frac{1}{2} + \frac{1}{2}\cos\frac{k\pi}{M}.  
\end{equation}
In such a case, the coefficients $\gamma^{[t_k]}_i$ 
and $\hat{\gamma}_l$ in (\ref{E:TchExpan}) are given by
\begin{equation}
\label{E:gammas}
\begin{array}{ll}
    \displaystyle \gamma^{[t_k]}_i = \frac{2-\delta_{i,M_{k}}}{M_{k}}\sum_{j=0}^{M_{k}}{''\,}
    \omega^{*}(s_j,t_k)^{-1}\cos\frac{ij\pi}{M_{k}},&\quad 0\leq i\leq M_{k},\\[3ex]
  \displaystyle \hat{\gamma}_l = \frac{2-\delta_{l,M}}{M}\sum_{k=0}^{M}{''\,}
    J(c,d;S_{M_{k}}) \cos\frac{lk\pi}{M},&\quad 0\leq l\leq M,
\end{array}
\end{equation}
where $\delta_{j,k}$ is the Kronecker delta, the double prime means that the first
and the last term of the sum are to be halved. The sets of coefficients (\ref{E:gammas})
can be very efficiently computed by means of the FFT with only $O\big(M_{k}\log(M_{k})\big)$
and $O\big(M\log(M)\big)$ arithmetic operations (cf.~\cite{Gen72} or \cite[Section 5.1]{DB08};
the authors recall that the FFT is not only fast, but also resistant to round-off errors).
The presented approach is very convenient from the practical point of view because
if the accuracy of the approximation (\ref{E:TchExpan}) is not satisfactory, then we may
double the value of $M_{k}$ (or $M$) and reuse the previously
computed results. The expansions (\ref{E:TchExpan}) are accepted if
\begin{equation}
\label{E:stop}
  \frac{\sum\limits_{i=M_{k}-3}^{M_{k}}|\gamma^{[t_k]}_i|}
       {\max\big\{1,\,\max\limits_{\,0\leq i\leq 3}|\gamma^{[t_k]}_i|\big\}} \leq 
         16\varepsilon
    \quad\, \mathrm{and} \quad\,
  \frac{\sum\limits_{i=M-3}^{M}|\hat{\gamma}_i|}
       {\max\big\{1,\,\max\limits_{\,0\leq i\leq 3}|\hat{\gamma}_i|\big\}} \leq
         256\varepsilon ,
\end{equation}
where $\varepsilon$ is the computation precision.

Here is the complete algorithm for efficient approximation of
the whole set of integrals $I_{\sbl j}$ for $\bl j\in\Omega_{n+m}^{\sbl c}$.
The functions (parameters) $a$, $b$, $c$, and $d$ are defined in (\ref{E:abcd}).
\begin{alg}[Numerical computation of the set of integrals
  $I_{\sbl j}$, $\bl j\in\Omega^{\sbl c}_{n+m}$]
\label{A:I-comp}
\
\begin{itemize}
\setlength{\itemindent}{-3.3ex}
\itemsep2pt
\item[] Let $M := M^*$, where $M^*$ is an arbitrary integer greater than 7.
\item[] \textbf{Phase I}. For $k\in\{0,1,\dots,M\}$ do the following Steps 1--6:\\[-2ex]
\begin{itemize}
\setlength{\itemindent}{5ex}
\setlength\itemsep{0.1ex}
\item[Step 1.] Compute $t_k$ according to \eqref{E:sjtk},
               and compute $w_i(t_k)$ in \eqref{E:omegastar} for $i\in\{0,1,\dots,n\}$.
\item[Step 2.] Let $M_{k} := M_{k}^{*}$, where $M_{k}^{*}$ is an arbitrary integer greater than 7.
\item[Step 3.] Compute the values $\hts\omega^{*}(s_j,t_k)^{-1}$ for $j\in\{0,1,\dots,M_{k}\}$,
               where $s_j$ is given by (\ref{E:sjtk}).
\item[Step 4.] Using the FFT, compute the coefficients $\gamma_i^{[t_k]}$ ($\,0\leq i\leq M_{k}$)
               defined in (\ref{E:gammas}).
\item[Step 5.] If the first condition of (\ref{E:stop}) is not satisfied, then set $M_{k} := 2 M_{k}$,
               compute the additional values $\hts\omega^{*}(s_j,t_k)^{-1}$ for $j\in\{1,3,5,\dots,M_{k}-1\}$,
               and go to Step 4.
\item[Step 6.] Compute the set of quantities
	\( W[t_k,j_1] := J\left(c(j_1),d(j_1);S_{M_{k}}\right),
               \) 
               by applying Algorithm \ref{A:PK_1d}, for
               $j_1\in\{c_1,c_1+1,\dots,N-c_2-c_3\}$, where $ N=n+m$.
\end{itemize}\vspace*{1ex}
\item[] \textbf{Phase II}. For $j_1\in\{c_1,c_1+1,\dots,N-c_3-c_2\}$ perform the following Steps 7--9:\\[-2ex]
\begin{itemize}
\setlength{\itemindent}{5ex}
\setlength\itemsep{0.1ex}
\item[Step 7.] Compute the coefficients $\hat{\gamma}_l$ ($\,0\leq l\leq M$)
               defined in (\ref{E:gammas}), by means of the FFT, using the stored values
               $W[t_k,j_1]$, $0\leq k\leq M$, in place of $J\big(c(j_1),d(j_1);S_{M_{k}}\big)$.
\item[Step 8.] If the second condition of (\ref{E:stop}) is not satisfied, then set
               $M := 2 M$, and repeat Steps 1--6 for $k\in\{1,3,5,\dots,M-1\}$.
\item[Step 9.] For $j_2\in\{c_2,c_2+1,\dots,N-c_3-j_1\}$,  compute the integrals
               \[
                 I_{\sbl j}\equiv I_{(j_1,j_2)} := A_{\sbla}\binom{N}{\bl j}\,
                   J\Big(a(\bl j),b(j_2);\hat{S}_{M}\Big)
               \]  
               using Algorithm \ref{A:PK_1d}.
\end{itemize}
\end{itemize}
\noindent \textbf{Output}: Set of the integrals $ I_{\sbl j}$
for $\bl  j\in\Omega^{\sbl  c}_{n+m}$.
\end{alg}
\begin{rem}
In Steps 4 and 7 of the above algorithm the coefficients $\gamma_i^{[t_k]}$
($\,0\leq i\leq M_{k}$) or $\hat{\gamma}_l$ ($\,0\leq l\leq M$) are recalculated
each time the value of $M_k$ or $M$ is doubled. Such a procedure is advised if we
use a system (like, e.g., Maple or Matlab) equipped with a fast built-in FFT
subroutine. If we are to program the FFT summation algorithm by ourselves,
it should rather be done in such a way that practically all results computed for
a previous value of $M_k$ or $M$ are reused
(cf., e.g.,\ \cite{Gen72}).
\end{rem}
In Table \ref{Tab:IntTime} we present the results of the efficiency test, where the
proposed quadrature (implemented in Maple) is compared to the Maple built-in integration
subroutine. We have used the B\'ezier surface form Example \ref{SS:example1} ($n=6$), and
set the parameters $m$ and $\bl c$ to several different values, to obtain collections
of integrals of different sizes (equal to $|\Omega_{n+m}^{\sbl c}|$). The experiment
was performed in the 64-bit version of Maple 16 on the computer equipped with
the $3.7$GHz i7 processor. All parameters $\alpha_i$ in (\ref{E:w}) were
set to $0$ (the efficiency of the proposed method does not depend on
$\bl\alpha$, but the Maple built-in integration subroutine
works most efficiently with this selection).

\begin{table}[!ht]
\newcolumntype{C}[1]{>{\centering\arraybackslash}m{#1}}
\newcommand{\hpo}{\hphantom{0}}
\begin{center}\vspace{-1.5ex}
	\caption{\small
Comparison of the computation times of the Maple library function
and the proposed adaptive quadrature (Algorithm \ref{A:I-comp}) in the
case of several collections of integrals (\ref{E:I}). The number of integrals
which are to be computed equals $|\Omega_{n+m}^{\sbl c}|$.}
								\label{Tab:IntTime}
\renewcommand{\arraystretch}{1.025}
\setlength\tabcolsep{1.5ex}
\vspace{0.4ex}
\begin{tabular}{C{14ex}C{24ex}C{24ex}}\hline
\vspace*{1.5ex}$\big|\Omega_{n+m}^{\sbl c}\big|$\vspace*{-1.75ex}
  & \multicolumn{2}{c}{computation time (in seconds)}\\\cline{2-3}
  & Maple library function & the proposed method\\\hline
  $\hpo\hpo1$ & $\hpo0.064$    & $0.30$ \\
  $\hpo\hpo3$ & $\hpo0.19\hpo$ & $0.30$ \\
     $\hpo10$ & $\hpo0.64\hpo$ & $0.32$ \\
     $\hpo28$ & $\hpo1.75\hpo$ & $0.37$ \\
     $\hpo91$ & $\hpo6.34\hpo$ & $0.43$ \\
        $276$ & $22.9\hpo\hpo$ & $0.59$ \\
        $990$ & FAILURE        & $0.89$ \\
\hline
\end{tabular}
\vspace{-2.5ex}
\end{center}
\end{table}

We have to keep in mind that Maple is an interpretative programming language
with a pretty slow code interpreter. Therefore, the $4.7$ times longer computation
time of our quadrature, compared to the computation time of the Maple library function,
in the case of $1$-element collection of integrals is in fact an excellent result.
The last collection of $990$ integrals ($n+m=42$) was too difficult to be
computed by the Maple built-in subroutine (in $14$-decimal digit
arithmetic, assumed during this test).


\section*{Appendix B: Hahn orthogonal polynomials}
						\label{S:AppB}
\setcounter{equation}0
\setcounter{subsection}0
\setcounter{thm}0
\renewcommand{\thethm}{B.\arabic{thm}}
\renewcommand{\thesection}{\Arabic{section}}
\renewcommand{\thesubsection}{B.\arabic{subsection}}
\renewcommand{\theequation}{B.\arabic{equation}}


The notation
\[
 \hyper rs {a_1,\ldots, a_r}{b_1,\ldots,b_s}z
 :=\sum_{k=0}^{\infty}\frac{(a_1)_k\cdots (a_r)_k}
                           {k!(b_1)_k\cdots (b_s)_k}\,z^k
\]
is used for the \textit{generalized hypergeometric series}
(see, e.g., \cite[\S2.1]{AAR99});
here $r,\,s\in\Zplus$, $z$, $a_i, b_j\in \C$,
and $(c)_k$ is the shifted factorial.
The \textit{Hahn polynomials} (see, e.g., \cite[\S1.5]{KS98})
\begin{equation}
	\label{E:Hahn1}
	h_l(t)\equiv h_l(t;a,b,M):=(a+1)_l(-M)_l\hyper32{-l,l+a+b+1,-t}{a+1,-M}1, 
\end{equation}  
where $l=0,1,\ldots,M$, $a,\,b>-1$, and $M\in\N$,  satisfy the  recurrence relation 
\begin{equation}
		\label{E:Hahn1-rec}
		h_{l+1}(t)=A_l(t,M)\,h_l(t)+B_l(M)\,h_{l-1}(t),
		\qquad l\ge0;\;	h_{0}(t)\equiv1;	\;h_{-1}(t)\equiv0,
	 \end{equation}
 with the coefficients
\begin{equation}	\label{E:Hahn1-rec-coeffs}
	A_l(t,M):=C_l\,(2l+s-1)_2\,t-D_l-E_l,\qquad 	
 B_l(M):=-D_l\,E_{l-1},		  	
\end{equation}
where $s:=a+b+1$,
$C_l:=(2l+s+1)/[(l+s)(2l+s-1)]$, $D_l:=C_l\,l(l+M+s)(l+b)$, and
	$E_l:=(l+a+1)(M-l)$.
\begin{rem}\label{R:Clenshaw}
	A linear combination of Hahn polynomials,
	\(
		s_N(t):=\sum_{i=0}^{N}\gamma_i\,h_i(t;a,b,M),
	\)	
	can be summed using the following \textit{Clenshaw's algorithm} (see, e.g., \cite[Thm 3.2.11]{DB08}).				
	Compute the sequence $V_0,V_1,\ldots,V_{n+2}$ from						
			$V_i:=\gamma_i+A_i(t;M)V_{i+1}+B_{i+1}(M)V_{i+2}$,  $i=N,N-1,\ldots,0$,
	with $V_{N+1}=V_{N+2}=0$, where the coefficients $A_i(t;M)$ and $B_i(M)$ are defined by \eqref{E:Hahn1-rec-coeffs}. Then $s_N(t)=V_0$.			 				 
\end{rem}



\small
\vspace{-4ex}
\begin{thebibliography}{99}
\itemsep0.45pt 
	
\bibitem{Ahl79}  L.V.~Ahlfors, Complex Analysis, 3rd Ed., McGraw-Hill, 1979.	

\bibitem{AAR99} G.E. Andrews, R. Askey, R. Roy,
                Special Functions,
                Cambridge Univ. Press, Cambridge,  1999.
                                                         
\bibitem{CW11}   J. Chen, G.J. Wang, Progressive-iterative approximation 
		for triangular B\'ezier surfaces,  Comp. Aided-Design 43 (2011) 889--895. 

\bibitem{DB08} G.~Dahlquist, A.~Bj\"orck, Numerical Methods in Scientific Computing,
  			vol. I, SIAM, Philadelphia, 2008.

\bibitem{Far86} G. Farin,  Triangular Bernstein-B\'ezier patches,
                Comput. Aided Geom. Design 3 (1986) 83--127.  
		
\bibitem{Far02} G. Farin, 
                Curves and Surfaces for Computer-Aided
                Geometric Design. A Practical Guide, 5th ed.,
                Academic Press,  Boston, 2002.
                
\bibitem{Gen72} W.M.~Gentleman, Implementing {C}lenshaw-{C}urtis quadrature --
                {II}. {C}omputing the cosine transformation,
                Comm. ACM {15} (1972) 343--346.

\bibitem{Hu13} Q.Q. Hu, An iterative algorithm for polynomial approximation of rational triangular
	       B\'ezier surfaces, Appl. Math. Comput. 219 (2013) 9308--9316.	       
	       
\bibitem{Kel07} P. Keller,  A method for indefinite integration of oscillatory and singular
		functions, Numer. Algor. {46} (2007)  219--251.	       
		

\bibitem{KS98} R. Koekoek, R.F. Swarttouw, 
               The Askey scheme of hypergeometric
               orthogonal polynomials and its $q$-analogue,  Rep. 98-17,
               Fac. Techn. Math. Informatics, Delft Univ. of Technology,
               Delft,  1998.
 
\bibitem{LKW15}  S. Lewanowicz,    P. Keller, P. Wo\'zny,
               B\'ezier form of dual bivariate Bernstein polynomials, 
               \texttt{arxiv:1510.08246 [math.NA]}  (2015).
       
\bibitem{LW06} S. Lewanowicz, P. Wo\'zny, 
	        Connections between two-variable
                Bernstein and Jacobi polynomials on the triangle, 
	        J. Comput. Appl. Math. 197 (2006) 520--533.               
               
\bibitem{LWK15} S. Lewanowicz, P. Wo\'zny,   P. Keller,
		Weighted polynomial  approximation of rational B\'ezier curves,
		\texttt{arXiv:1502.07877 [math.NA]} (2015).	    	        

\bibitem{LWK12} S. Lewanowicz, P. Wo\'zny,   P. Keller,
                Polynomial  approximation of rational B\'ezier curves with constraints,
                Numer. Algor. 59 (2012) 607--622.	        
	       
\bibitem{Rab05} A. Rababah,  Distances with rational triangular B\'ezier surfaces,
	        Appl. Math. Comp. 160 (2005)  379--386.

\bibitem{Riv90} T.J.~Rivlin, Chebyshev Polynomials:
  			From Approximation Theory to Algebra and Number Theory,
  			2nd ed., Wiley, New York, 1990.
	        
\bibitem{Sha13}	R. Sharma,   Conversion of a rational polynomial in triangular
		Bernstein-B\'ezier form to polynomial in triangular  Bernstein-B\'ezier form,  
		Internat. J. Comput. Appl. Math. 8 (2013) 45--52.
    		              
\bibitem{WL10}  P. Wo\'zny, S.  Lewanowicz,   
		Constrained multi-degree reduction of triangular B\'ezier
		surfaces using dual Bernstein polynomials,
                J. Comput. Appl. Math. 235 (2010) 785--804.                
		
\bibitem{WL08}  P. Wo\'zny, S.  Lewanowicz,   
		Multi-degree reduction of B\'ezier
		curves with constraints, using dual Bernstein basis
		polynomials, Comput. Aided Geom. Design 26 (2009) 566--579.

				
\bibitem{XW09} H.X. Xu, G.J. Wang, Approximating rational  triangular B\'ezier 
		surfaces by polynomial triangular B\'ezier surfaces,  		 
		J. Comput. Appl. Math. 228 (2009) 287--295.
		
\bibitem{XW10} H.X. Xu, G.J. Wang,	New algorithm of polynomial triangular B-B surfaces
		approximation to rational triangular B-B surfaces, J. Inform. Comput. Sci. 7 (2010) 725--738.
		
\bibitem{ZW06} L. Zhang, G.J. Wang, An effective algorithm to approximate rational
		triangular B-B surfaces using polynomial forms, Chinese J. Computers 29 (2006) 2151--2162
		(\textit{in Chinese}).		
		
\end{thebibliography}
\end{document}